\documentclass{amsart}
\usepackage{amssymb, amsmath}

\def\largerightarrow{   -\negthinspace\negthinspace -\negthinspace\negthinspace
                                       \negthinspace\longrightarrow }

\begin{document}
\Large
\title{}
\vbox{\hfil {\Large\bf SUPER OPERATOR SYSTEMS, STRONG NORMS }\hfil}
\vbox{\hfil {\Large\bf AND OPERATOR TENSOR PRODUCTS }\hfil}
\author{U. Haag}
\date{\today}
\maketitle
\noindent
\begin{abstract}
A notion of super operator system is defined which generalizes the 
usual notion of operator systems to include certain unital involutive operator
spaces which cannot be represented completely isometric as a concrete 
operator system on some Hilbert space. They can nevertheless be represented 
by bounded operators on a standard $\mathbb Z_2$-graded Hilbert space 
equipped with a superinvolution. We apply this theory to investigate on the relation 
between certain tensor products defined for operator spaces and $C^*$-algebras, 
such as the projective tensor product, the Haagerup tensor product and the maximal 
$C^*$-tensor product.
\end{abstract}
\par\bigskip\bigskip\bigskip\bigskip\bigskip\noindent
{\bf 0. Introduction}
\par\bigskip\noindent
As the title indicates this paper is primarily concerned with the description and examination of certain operator spaces with some additional structure we call {\it super operator systems}. If an operator space 
can be regarded as the quantization of a Banach space, then a unital operator space is sort of an operator space with a backbone. A super operator system is a unital operator space which posesses an antilinear involution such that for each natural number $\, n\, $ its composition with the transposition operation on $n\times n$-matrices is isometric on the $n$-fold amplification of the operator space (tensor product with $\, M_n ( \mathbb C ) $). Also the involution is required to fix the identity element. Because the term {\it super} is sometimes overused we feel obliged to give some justification. 
Following the general abstract characterization of an operator space due to Ruan (cf. \cite{E-R})  abstract characterizations of a unital operator space have been provided in \cite{B-N} and \cite{Bl-Ne}. 
One knows that given any concretely represented unital operator space, its enveloping operator system is uniquely determined up to complete $*$-isometry regardless of the particular representation so that it is really determined by the abstract structure of a unital operator space. Therefore a super operator system has a unique enveloping operator system which is in fact ${\mathbb Z}_2$-graded. 
The (completely isometric $*$-linear) grading comes in by composing the antilinear involution with the antilinear adjoint operation of the enveloping operator system. It is then easy to see that any abstract 
super operator system has a concrete realization on ${\mathbb Z}_2$-graded 
Hilbert space which is $*$-linear with respect to the product of the grading automorphism plus the adjoint operation which is called the {\it super involution}. On the other hand it is not always possible to get a $*$-linear representation on ungraded Hilbert space (i.e. if and only if the super operator system is in fact an operator system). Since there are interesting super operator systems appearing quite naturally, for example by considering certain operator tensor products of operator systems or $C^*$-algebras, which are not operator systems by themselves one might say that ${\mathbb Z}_2$-gradings occur naturally in mathematics 
and especially in operator theory. The term super operator system then appears quite appropriate for this generalization of the concept of an operator system.  One is also reminded of certain concepts in 
mathematical physics, especially the notion of supersymmetry (c.f. the article  \cite{H-L-S} of R. Haag, J. Lopuszanski and M. Sohnius). The author knows to little about these matters to decide whether there might be some deeper applications of the theory of super operator systems in physics but hopes that physicists will feel challenged to consider this question in more detail. This paper is self-contained to a large extent, using only standard results from the theory of operator spaces which can all be found together with further references to original papers in the books \cite{E-R} by E. Effros and Z.-J. Ruan and 
\cite{Pa} by V. Paulsen,  
so we have chosen these two as a general reference frame for this text. The article is organized into three sections: the first gives the basic definitions and results for super operator systems and super 
$C^*$-algebras, including the notion of superpositivity, the second section is an excursion to (super) operator spaces equipped with what is called a strong matrix norm in our terminology, it does not have any direct applications in the sequel so the 
reader mainly interested in the results of the last section may well skip the second one, but the feeling of the author is that some of the results (and primarily the so called Strong Lemma) might have very interesting applications in the future yet to be discovered (as was intended in the writing) so we choose to maintain 
the section as it is leaving the judgement whether or not it is worthwile reading to the recipient. 
Finally, the third section may be considered as the main part since it gives some of the more profound applications of the theory, most of them dealing with operator tensor products of various kind, and centered around the Super-Christensen-Effros-Sinclair Theorem (SCES-Theorem for short) which gives some sufficient conditions for the action of a $C^*$-algebra or topological group on an operator system to be unitarily implementable (strongly completely contractive). 
\par\bigskip\noindent 
{\bf 1. Super operator systems and superpositivity}
\par\bigskip\noindent 
Let $\, {\mathcal H}\, $ denote a 
Hilbert space of a given finite or infinite dimension, usually we assume that $\, {\mathcal H}\, $ is infinite but countably generated. The standard $\, {\mathbb Z}_2$-graded Hilbert space is the direct sum 
of two copies of $\, {\mathcal H}\, $ one of which is supposed to have even and the other odd grading, 
$\, \widehat {\mathcal H}\, = {\mathcal  H}^{even} \oplus {\mathcal H}^{odd}\, $. Thus the grading operator 
on $\,\widehat {\mathcal H }\, $ is given by the matrix 
$$ \epsilon =
 \begin{pmatrix}
1 & \>\> 0 \\
0 & -1 
\end{pmatrix} $$ 
with respect to this decomposition. Conjugating bounded operators on $\,\widehat {\mathcal H}\, $ by 
$\,\epsilon\, $ defines an order two automorphism compatible with the adjoint operation. One then defines the superinvolution on $\, \mathcal B ( \widehat {\mathcal H } )\, $ as product of the grading automorphism and the adjoint operation. It is clear that the superinvolution is an isometric antilinear antiautomorphism of $\, \mathcal B ( \widehat {\mathcal H } )\, $. In the following the notation 
$\, x \mapsto x^*\, $ will always allude to the superinvolution unless stated otherwise. To prevent mixing up the super involution with the ordinary adjoint operation the latter is denoted 
$\, x\mapsto \overline x\, $. By definition, 
a {\it super-$C^*$-algebra} consists of a norm-closed subalgebra of 
$\, \mathcal B ( \widehat {\mathcal H })\, $ which is invariant under the superinvolution. Note that a 
super-$C^*$-algebra need not be $\, {\mathbb Z}_2$-graded by itself, i.e. it is not required to be closed under the grading automorphism (and as a consequence also not by the ordinary adjoint operation). It is a certain type of operator algebra which in special cases reduces to an ordinary 
$C^*$-algebra (possibly with ${\mathbb Z}_2$-grading).  Much in the same way one defines a (concrete) {\it super operator system} to consist of a norm-closed subspace $\,\mathfrak X\, $ of 
$\, \mathcal B ( \widehat {\mathcal H} )\, $ which contains the unit element and is invariant under the superinvolution. Such a subspace is naturally equipped with the induced operator matrix norms on the 
collection of subspaces $\, M_n ( \mathfrak X ) \subseteq M_n ( \mathcal B ( \widehat {\mathcal H} ))\, $.
$\, \mathfrak X\, $ may or may not contain the grading operator $\, \epsilon\, $. If it does we will refer to 
$\,\mathfrak X\, $ as a {\it complete super operator system}. Let $\, A\, $ be a (for the time being concrete) 
super $C^*$-algebra. We want our representation to be as simple as possible. Given any $*$-representation of $\, A\, $ on a 
$\ \mathbb Z_2$-graded Hilbert space, consider the enveloping $\,\mathbb Z_2$-graded $C^*$-algebra 
$\,\widehat A\, $ and let $\,\check A\, $ denote its underlying trivially graded $C^*$-algebra. Any (faithful) graded representation of $\,\widehat A\, $ entails a (completely) isometric $*$-representation of $\, A\, $. Now assume given a representation of $\,\check A\, $ on ordinary ungraded Hilbert space. Let $\, C_1\, $ denote the complex Clifford algebra on one odd generator. One has a 
natural $\,\mathbb Z_2$-graded isomorphism $\, \check A \otimes C_1 \simeq \widehat A \otimes C_1\, $ so that the natural embeddings $\, \widehat A \hookrightarrow \widehat A \otimes C_1\, $  and 
$\, \check A \otimes C_1 \hookrightarrow \widehat M_2 ( \check A )\, $ ($\, {\widehat M}_2 ( \check A )\, $ is the $C^*$-algebra $\, M_2 ( \check A )\, $ but with standard off-diagonal grading) 
determine a $\,\mathbb Z_2$-graded representation of $\,\widehat A\, $ via the ungraded representation of $\,\check A\, $. Then any hermitian (= superselfadjoint) element of $\,\widehat A\, $ takes the form
$$ \begin{pmatrix} 
a & i b \\
i b & a \end{pmatrix} $$
with $\, a\, $ and $\, b\, $ selfadjoint. Let us call such a representation  {\it of standard form}.
Given a contractive unital $*$-homomorphism of one unital super-$C^*$-algebra 
$\, A\, $ into a unital super-$C^*$-algebra $\, B\, $, both represented in standard form, one obtains an induced contractive $*$-homomorphism from $\, A \otimes C_1\, $ to $\, B \otimes C_1\, $ where the tensor product is endowed with the obvious multiplication and product type superinvolution. Actually, 
$\, A \otimes C_1\, $ is represented on the same $\,\mathbb Z_2$-graded Hilbert space, sending the degree one  (ordinary) self-adjoint generator $\, F\, $ of $\, C_1\, $ with  $\, F^2 = 1\, $ to the matrix 
$$ \begin{pmatrix} 0 & 1 \\ 1 & 0 \end{pmatrix}  $$
owing to the fact that $\, A\, $ is in standard form.
Consider the super-$C^*$-algebras $\, \alpha ( A )\, $ and $\, \alpha ( B  )\, $ (we will always let 
$\,\alpha\, $ denote the grading automorphism on the ambient algebra of bounded operators). Since 
$\, \alpha\, $ is an isometry one obtains a contractive $*$-homomorphism from
$\, \alpha ( A  )\, $ to $\, \alpha ( B  )\, $. We would like to glue these two $*$-homomorphisms together 
over the intersection $\, ( A \cap \alpha ( A ) )\, $ and extend to a $\,\mathbb Z_2$-graded 
$*$-homomorphism of the enveloping $\,\mathbb Z_2$-graded $C^*$-algebras 
$\, \widehat A\, $, $\,\widehat B\, $. In order to do this one must ensure that both homomorphisms agree on the intersection which is clearly a 
$\,\mathbb Z_2$-graded $C^*$-algebra by itself being closed under superinvolution as well as the grading automorphism. A necessary condition is that the restricted kernel of $\, \varphi\, $ is invariant under the grading automorphism. This turns out to be automatic from the fact that the restricted kernel is necessarily a $*$-invariant ideal and the Lemma below.  Then the image of the intersection as above under a contractive $*$-homomorphism is isometrically isomorphic with a $\,\mathbb Z_2$-graded $C^*$-algebra. 
A contractive $*$-representation of such an algebra is also selfadjoint with respect to the ordinary involution and hence graded by the following Lemma (this result is well known and included only for convenience).
\par\bigskip\noindent
{\bf Lemma 1.}\quad Let $\, \varphi : \widehat A \rightarrowtail \mathcal B ( \widehat {\mathcal H} )\, $ be a contractive unital $*$-representation (with respect to superinvolution) of a $\,\mathbb Z_2$-graded unital $C^*$-algebra. Then 
$\, \varphi\, $ is a graded $C^*$-representation.  
\par\bigskip\noindent
{\it Proof.}\quad Any unital contractive representation of an ordinary $C^*$-algebra is selfadjoint. This follows for example from the fact that an operator in $\,\mathcal B ( \mathcal H )\, $ is selfadjoint of norm less or equal to one if and only if  
$\, \Vert x - i t {\bf 1} \Vert \leq \sqrt{1 + t^2}\, $ for all $\, t \in \mathbb R\, $ (c.f. \cite{E-R}, Lemma 
A.4.2) \qed
\par\bigskip\noindent
Consider the $C^*$-algebra $\, C^* ( \bf 1 , \epsilon ) \simeq \mathbb C \oplus \mathbb C\, $ generated by the unit and the grading operator and let $\,\omega\, $ be any unitary operator of this algebra. One can define an isometric grading preserving involution on $\,\mathcal B ( \widehat{\mathcal H} )\, $ by 
the formula $\, x \mapsto x^{\omega } := \omega\cdot \overline x \cdot\omega\, $ where $\, x \mapsto \overline x\, $ denotes the ordinary adjoint operation. In particular $\, \overline x = x^{\bf 1}\, $ and 
$\, x^* = x^{\epsilonÊ}\, $. An element such that $\, x^{\omega } = x\, $ will be called 
$\omega$-hermitian.  
There is an analogue of the cited result for $\omega$-hermitian elements, namely
\par\bigskip\noindent
{\bf Lemma 2.}\quad Let $\, x\, $ be an operator in $\, \mathcal B ( \widehat {\mathcal H } )\, $ with specified grading operator $\,\epsilon\, $ as above. Then $\, x\, $ is $\omega$-hermitian of norm less or equal to one if and only if for every $\, t \in \mathbb R\, $ the following inequality holds
$$ \Vert x - i t \omega \Vert \leq \sqrt{ 1 + t^2} \> . $$
\par\bigskip\noindent
{\it Proof.}\quad Let $\, x\, $ be $\omega$-hermitian of norm less or equal than one. Choose a unitary square root $\,\gamma\, $ of $\, {\omega }^{-1}\, $, i.e. $\, {\gamma }^2 = {\omegaÊ}^{-1}\, $. With respect to the matrix decomposition given by the grading one has 
$$ \omega = \begin{pmatrix} {\omega }_0 & 0 \\ 0 & {\omegaÊ}_1 \end{pmatrix}\quad , \quad 
\gamma = \begin{pmatrix} {\gamma }_0 & 0 \\ 0 & {\gamma }_1 \end{pmatrix}\quad ,\quad 
x = \begin{pmatrix} a & b \\ c & d \end{pmatrix}\> . $$
The condition that $\, x\, $ is $\omega$-hermitian means that $\, a = {\omega }_0^2\, a^*\, ,\, 
b = {\omega }_0\, {\omega }_1\, c^*\, ,\, d = {\omega }_1^2\, d^*\, $. 
Multiplying the inequality $\, \Vert x - i t \omega \Vert \leq \sqrt{ 1 + t^2}\, $ by $\, \gamma\, $ from the left and from the right one obtains 
$$ \Vert \widetilde x - i t {\bf 1} \Vert \leq \sqrt{ 1 + t^2 } $$
with $\, \widetilde x = \gamma\cdot x\cdot \gamma\, $ an ordinary selfadjoint element whenever
$$ \Vert x \Vert \leq 1 $$
by Lemma A.4.2 of \cite{E-R}.  The reverse direction is much the same\qed                                                                                            
\par\bigskip\noindent
Note however that the $\omega$-involution is (anti)multiplicative and unit preserving only if $\,\omega\, $ is selfadjoint. In the applications we will usually have $\, {\omega }_0 = 1\, $.
In view of  Lemma 1 we may glue the two maps as above over the common intersection of their domains to obtain a $\,\mathbb Z_2$-graded $*$-linear map $\, \phi : ( A + \alpha ( A ) )  \rightarrow 
( B + \alpha ( B ) ) \, $ which restricted to either $\, A\, $ or $\, \alpha ( A )\, $ is a $*$-homomorphism. The extension is contractive by \cite{Pa}, Proposition 2.12.. 
Since $\,\phi\, $ is also unital, it is positive in the ordinary sense. The same construction carries over to matrices of arbitrary size whenever $\, \varphi \, $ is completely contractive. Then $\,\phi\, $ is completely positive. Note that we have not made explicit use of the fact that 
$\,\varphi\, $ is an algebra homomorphism. Restricting to the $*$-linear structure the same arguments yield for every $*$-invariant subspace $\,\mathfrak X\, $ of a super-$C^*$-algebra canonically an enveloping $\,\mathbb Z_2$-graded operator system 
$\,\widehat {\mathfrak X}\, $. In particular starting with a $\,\mathbb Z_2$-graded operator system, its structure as an operator system is uniquely determined by its structure as a super operator system. But the latter 
contrary to the former also determines the $\,\mathbb Z_2$-grading and is thus {\it finer} than the operator system structure which might be a bit of a surprise.  We would like to extend not only the $*$-linear, but also the multiplicative structure to the graded setting.  However the enveloping graded 
$C^*$-algebra need not be unique, but may depend on the particular representation. It is known that for any unital operator algebra $\,\mathcal A\, $ there is a minimal choice of enveloping $C^*$-algebra 
$\, C_{min}^* ( \mathcal A )\, $ which has the property that any other $C^*$-algebra generated by 
$\,\mathcal A\, $ in some particular representation canonically surjects onto 
$\,C_{min}^* ( \mathcal A )\, $ (Hamana's Theorem, compare \cite{Pa} Theorem 15.16). Correspondingly there is a maximal choice of enveloping 
$C^*$-algebra which is obtained by taking the direct sum of all (spatial equivalence classes of) 
completely contractive unital representations of $\,\mathcal A\, $. Denote this algebra by 
$\, C_{max}^* ( \mathcal A )\, $. This construction has the advantage that it is functorial, i.e. 
given any completely contractive homomorphism $\,\mathcal A \rightarrow \mathcal B\, $ one gets an induced $*$-homomorphism $\, C_{max}^* ( \mathcal A ) \rightarrow C_{max}^* ( \mathcal B )\, $. The following simple example illustrates that $\, C_{max} ( \mathcal A )\, $ might be different from 
$\, C_{min}^*( \mathcal A )\, $.
\par\bigskip\noindent
{\it Example.}\quad Let 
$$ x\> =\> \begin{pmatrix} \>\>\> 1 & \>\>\> 1 \\ -1 & -1 \end{pmatrix} \in \widehat M_2 ( \mathbb C ) $$
and $\, A_0\, $ be the commutative super-$C^*$-algebra generated by $\, x\, $ and $\,\bf 1\, $. Then 
the minimal enveloping graded $C^*$-algebra is clearly equal to $\,\widehat M_2 ( \mathbb C )\, $. On the other hand, since $\, x\, $ is nilpotent, there exists a nontrivial $*$-homomorphism 
$\, A_0 \rightarrow \mathbb C\, $ that cannot be extended to $\,\widehat M_2 ( \mathbb C )\, $ which is simple, so that $\, \widehat A_{0 , max} = C_{max}^* ( A_0 ) \, $ must be an extension of 
$\,\widehat M_2 ( \mathbb C )\, $.
\par\bigskip\noindent
Let us return to the concept of super operator system for which one would like to have some abstract characterization. There is some hidden structure which so far has remained unrevealed. This is connected with the concept of {\it superpositivity} and of {\it $\epsilon$-positivity}. As we have seen one recaptures from a given super operator system $\,\mathfrak X\, $ a $\,\mathbb Z_2$-graded operator system $\,\widehat {\mathfrak X}\, $ which entails the usual notion of matrix ordering by {\it positive} elements. In a way the grading operator $\,\epsilon\, $ can be viewed as a sort of second order unit with respect to an ordering by so-called $\epsilon$-positive elements. Define the {\it graded spectrum} of an element 
$\, x\in \mathfrak X\, $ to be the subset $\,\widehat {\sigma} ( x ) \subseteq \mathbb C\, $ such that 
$\,\lambda \in \widehat {\sigma} ( x )\, $ iff the operator $\, x - \lambda\epsilon\, $ is not invertible in 
$\,\mathcal B ( \widehat {\mathcal H} )\, $. Of course this definition depends on the particular representation. For example it is clear that the graded spectrum of any element in a standard representation of $\,\mathfrak X\, $ is invariant with respect to the operation 
$\,\lambda \mapsto -\lambda\, $ whereas if $\,x\, $ is a positive element which is homogenous of even degree, there exists a representation such that 
$\, \widehat {\sigma} ( x )\, $ coincides with the ordinary spectrum in that representation which is purely positive. There are however some invariants. For example it follows from the proof of Lemma 2 that multiplication with the matrix 
$$ \kappa\> =\> \begin{pmatrix} 1 & 0 \\ 0 & i \end{pmatrix} $$ 
from the left and right defines a complete
isometry of a complete super operator system $\,  \mathfrak X \, $ containing the grading operator with an ordinary operator system $\, \check {\mathfrak X}\, $ such that the grading operator $\,\epsilon\, $ transforms into the unit element. Thus the norm of any hermitian element is given by the spectral radius of its graded spectrum. For a given complete super operator system $\,\mathfrak X\, $ define a matrix order on the hermitian elements by the (convex Archimedean) cones $\, {\mathcal C}^{\epsilonÊ}_{n , +} \subseteq M_n ( \mathfrak X )\, $ of elements with positive graded spectrum (= $\epsilon$-positive elements). From the completely isometric  order isomorphism of 
$\,\mathfrak X\, $ with $\, \check {\mathfrak X}\, $ one deduces that any hermitian element is the difference of two $\epsilon$-positive elements, i.e. 
$\, {M_n ( \mathfrak X )}^h = {\mathcal C}^{\epsilon }_{n , +} - {\mathcal C}^{\epsilon }_{n , +}\, $, and 
$\, {\mathcal C}^{\epsilon }_{n , +} \cap - {\mathcal C}^{\epsilon }_{n , +} = \{ 0 \}\, $. In the light of this new perspective we have to rethink our definition of a super-$C^*$-algebra, which in any given representation can be completed to a complete super operator system by adjoining the grading operator. Of course if the algebra already contains the grading operator it is necessarily a $\,\mathbb Z_2$-graded $C^*$-algebra, but this requirement seems to restrictive. Instead, one may think of the grading operator as a  virtual graded unit which is not part of the multiplicative algebraic structure of 
$\, A\, $, but only of the $*$-linear 
structure with induced $\epsilon$-ordering on the underlying complete super operator system. This will be called a {\it complete super $C^*$-algebra} as opposed to an ordinary super $C^*$-algebra given by the $*$-invariant (unital) operator algebra $\, A\, $. The 
$\epsilon$-order is not intrinsic in the $*$-multiplicative structure of a given super-$C^*$-algebra, and has to be specified separately. We will now try to give a definition which is independent of a particular representation. An abstract operator space with specified element 
$\, \bf 1\, $ which can be represented on Hilbert space such that the specified element corresponds to the unit element of $\,\mathcal B ( \mathcal H )\, $ will be called a {\it unital operator space}. It follows from \cite{Pa}, Prop. 3.5. that the enveloping operator system of a unital operator space is well defined irrespective of the particular unital representation. Any (ordinary) selfadjoint linear subspace of 
$\,\mathcal B ( \mathcal H )\, $ which does not necessarily contain the unit element
will be called a {\it pre-operator system}. If one wants an abstract definition of a (complete) super operator system one should first extract the relevant components of the structure of a given concrete (complete) super operator system. One starts with a unital operator space which is selfadjoint with respect to some isometric involution and contains beside the unit another specified element $\,\epsilon\, $, the graded unit. The graded unit has several aspects. First of all, being the grading operator in $\, \mathcal B ( \widehat{ \mathcal H} )\, $, it is a selfadjoint unitary and determines 
the grading automorphism on the completed $\mathbb Z_2$-graded operator system 
$\,\widehat {\mathfrak X}\, $ which is obtained through conjugation with $\,\epsilon\, $. Now the property of being invertible is not necessarily inherent in the abstract definition of an operator system (the spectrum of an element might differ with the concrete representation). But, attached to the operator system $\,\widehat {\mathfrak X}\, $ in a canonical manner, is its injective envelope 
$\, I ( \widehat{\mathfrak X} )\, $ which is a unital $C^*$-algebra. If $\,\widehat{ \mathfrak X}\, $ is represented on a Hilbert space $\, \mathcal H\, $, then 
$\, \widehat {\mathfrak X} \subseteq I ( \widehat{ \mathfrak X} ) \subseteq \mathcal B (  \mathcal H )\, $ and there exists a completely positive projection $\, \varphi : \mathcal B ( \mathcal H ) \rightarrow 
\mathcal B ( \mathcal H )\, $ with range $\, I ( \widehat{\mathfrak X} )\, $ such that 
$\,{\varphi }^2 = \varphi\, $, and multiplication in $\, I ( \widehat{\mathfrak X} )\, $ is given by the formula 
$\, x * y = \varphi ( x\cdot y )\, $. In particular, since $\,\epsilon\cdot\epsilon = 1\, $ one has 
$\,\epsilon * \epsilon = 1\, $ so that $\,\epsilon\, $ is a selfadjoint unitary of even degree in 
$\, I ( \widehat {\mathfrak X} )\, $. The grading automorphism extends from $\,\widehat{\mathfrak X}\, $ to 
$\, I ( \widehat{\mathfrak X} )\, $ by letting $\, \alpha ( x ) = \varphi ( \overline{\alpha } ( x ) )\, $ where 
$\,\overline{\alpha}\, $ is any completely positive extension of $\,\alpha\, $ to 
$\,\mathcal B ( \mathcal H  )\, $ (which exists from Arvesons extension theorem); of course in a 
$\mathbb Z_2$-graded representation the extension is obvious. The point is that no matter what extension, it is uniquely determined on $\, I ( \widehat{ \mathfrak X} )\, $ from the fact that 
$\,Ê{\alpha }^2 = id\, $ on $\,\widehat {\mathfrak X}\, $ (compare \cite{Pa}, chap. 15). From this and the identity 
$\, \epsilon * x * \epsilon = \varphi ( \epsilon\cdot x\cdot\epsilon )\, $ one gets that $\,\alpha\, $ is implemented in the abstract $C^*$-algebra $\, I ( \widehat{\mathfrak X} )\, $ through conjugation with 
$\, \epsilon\, $. We gather these results in the following definition. 
\par\bigskip\noindent
{\bf Definition 1.}\quad (i) $\> $ An {\it abstract complete super operator system} is a 
quadruple $\, ( \mathfrak X\, ,\, *\, ,\, \bf 1\, ,\, \epsilon )\, $ where $\, \mathfrak X\, $ is completely isometric to an abstract operator system with order unit $\epsilon\, $ (the graded unit) and completely antiisometric involution $\, *\, $ (which means that the involution composed with the transposition of matrices is isometric for each $\, M_n ( \mathfrak X )\, $), determining a matrix order given by a sequence of Archimedean cones 
of {\it $\epsilon$-positive} elements 
$\,\{ {\mathcal C}^{\epsilon }_{n , +}Ê\} \subseteq M_n ( \mathfrak X )^h\, $, and containing another specified element $\,\bf 1\, $ invariant under the involution, the unit of $\,\mathfrak X\, $. It then follows that its injective envelope 
$\, I ( \mathfrak X , \epsilon )\, $ is a $C^*$-algebra with unit $\,\epsilon\, $ and uniquely determined multiplication.  It is assumed that $\,\bf 1\, $ is a selfadjoint unitary for this multiplicative structure. Any 
$*$-invariant subspace of a complete abstract super operator system containing the unit $\, \bf 1\, $ is referred to as an (ordinary) super operator system.
\par\noindent
A {\it supermorphism} of super operator systems is a (completely) bounded 
$*$-linear map $\, s: \mathfrak X \rightarrow \mathfrak Y\, $. $\, s\, $ is {\it $\epsilon$-unital}, resp. 
{\it unital}, iff it maps the graded unit, resp. unit, of $\,\mathfrak X\, $ to the graded unit, resp. unit, of 
$\,\mathfrak Y\, $. 
$\, s\, $ is {\it (completely) $\epsilon$-positive} iff 
$\, s\, $ ($\, s_n\, $) maps $\epsilon$-positive elements into $\epsilon$-positive elements where 
$\, s_n : M_n ( \mathfrak X ) \rightarrow M_n ( \mathfrak Y )\, $ is the induced map.
\par\smallskip\noindent
(ii) An {\it abstract  (unital) complete super-$C^*$-algebra} is a pair $\, ( A , \epsilonÊ )\, $ where 
$\, A\, $ is a 
(unital) operator algebra equipped with a completely antiisometric antimultiplicative involution, also referred to as an (ordinary) {\it abstract  super-$C^*$-algebra}  which together with the matrix norms extends to the linear space 
$\, \mathfrak A = A + \mathbb C\cdot\epsilon\, $ such that $\, \mathfrak A\, $ is a complete super operator system with graded unit $\,\epsilon\, $. A {\it complete superhomomorphism} 
$\,\phi : ( A , {\epsilonÊ}_A ) \rightarrow ( B , {\epsilon }_B )\, $ is  a (completely) bounded supermorphism from $\, \mathfrak A\, $ to $\,\mathfrak B\, $ which restricts to
a (completely) bounded $*$-homomorphism from $\, A\, $ to $\, B\, $. 
\par\bigskip\noindent
The case that $\,\epsilon\, $ is contained in $\, A\, $ is not excluded in the above definition of a complete  super-$C^*$-algebra, but not assumed either. In (i) the case that $\,\epsilon\, $ equals $\,\bf 1\, $ is also not excluded. In this case $\,\mathfrak X\, $ becomes an operator system in the usual sense. 
\par\bigskip\noindent
{\it Example.}\quad As an example one checks that any two isometric standard representations of a given super-$C^*$-algebra $\, A\, $ extend to an isometric supermorphism of $\epsilon $-hermitian parts of the completed concrete super-$C^*$-algebras. The task is, given two different standard $*$-representations $\,\rho\, ,\, \sigma\, $ of $\, A\, $ to extend the $*$-isometry 
$\, \rho ( A ) \rightarrow \sigma ( A )\, $ to a $*$-isometry 
$\,  \rho ( A )^h + \mathbb R\cdotÊ{\epsilon}_{\rho} \rightarrow \sigma ( A )^h + \mathbb R\cdot {\epsilon }_{\sigma }\, $. First of all, since both representations are standard the grading operators 
$\, {\epsilon }_{\rho}\, $ and $\, {\epsilon }_{\sigma }\, $ are not contained in the respective images of 
$\, A\, $. Since the graded spectrum of each hermitian element $\, x\, $ is invariant by the inverse map, one easily calculates the norm of $\, x + \lambda \epsilon\, $ for each $\, \lambda \in \mathbb R\, $ and finds that it is independent of the representation. Thus the natural extension of 
$\, \rho ( A )^h \rightarrow \sigma ( A )^h\, $ to the hermitian parts of the completed super operator systems is an isometry, and $\epsilon$-unital, hence (bi-)$\epsilon$-positive.
\par\bigskip\noindent
We next show that any abstract complete super operator system (and abstract super-$C^*$-algebra) can be represented as a concrete super operator system on $\mathbb Z_2$-graded Hilbert space (concrete super-$C^*$-algebra). It is clear from the Choi-Effros representation theorem , that the structure of operator system with respect to the graded unit 
$\,\epsilon\, $ can be represented, so the problem is to find a representation which also suits the structure of unital operator space with respect to the unit $\,\bf 1\, $. Thanks to the discussion preceding the definition the proof is very simple (the incomplete case is deferred until a better understanding of abstract super operator system has been reached).
\par\bigskip\noindent
{\bf Theorem.} (Complete super operator systems representation theorem)\quad 
(i) Any abstract complete super operator system 
$\, ( \mathfrak X\, ,\, *\, ,\, \bf 1\, ,\, \epsilon )\, $ has a faithful unital completely isometric $*$-representation on a $\mathbb Z_2$-graded Hilbert space $\,\widehat {\mathcal H}\, $ (with respect to superinvolution) such that $\,\epsilon\, $ corresponds to the grading operator. It then follows that the enveloping operator system of 
$\, ( \mathfrak X\, ,\, \bf 1 )\, $ is well defined and decomposes as a direct sum 
$\,\widehat{\mathfrak X} = {\mathfrak X}_0 \oplus {\mathfrak X}_1\, $, where
$\,{\mathfrak X}_0\, $ is the operator system  generated by the set 
$\, \{ x + \overline x\,\vert\, x = x^*,\, x\in \mathfrak X \}\, $, and the pre-operator system 
$\, {\mathfrak X}_1\, $ is generated by the set 
$\,\{ x - \overline x\,\vert\, x = x^*,\, x \in \mathfrak X \}\, $. 
Moreover,
the two natural projections  
$\, p_0 : \widehat {\mathfrak X} \twoheadrightarrow {\mathfrak X}_0\, $ and
$\, p_1 : \widehat {\mathfrak X} \twoheadrightarrow {\mathfrak X}_1\, $ are complete quotient mappings with respect to the norms determined by the natural inclusions into 
$\,\widehat {\mathfrak X}\, $.
\par\smallskip\noindent
(ii) Any abstract super-$C^*$-algebra has a completely isometric unital $*$-
representation on 
$\mathbb Z_2$-graded Hilbert space. 
\par\bigskip\noindent 
{\it Proof.}\quad Consider the $\mathbb Z_2$-graded $C^*$-algebra $\, I ( \mathfrak X , \epsilon )\, $, the grading automorphism being defined as conjugation with the image of $\,\bf 1\, $. It has a unital faithful $C^*$-representation $\,\pi\, $ on some Hilbert space $\,\widehat{\mathcal H}\, $ which upon regarding the image of 
$\,\bf 1\, $ as a grading operator inherits a $\,\mathbb Z_2$-grading. Then consider the element 
$\, \kappa = \pi \left( ( {\bf 1} + \epsilon ) / 2 + i ( {\bf 1} - \epsilon ) / 2Ê\right)\, $. The map 
$\, x \mapsto \kappa\cdot x\cdot\kappa\, $ of $\, \mathcal B ( \widehat{\mathcal H} )\, $ exchanges the images of $\, \bf 1\, $ and $\,\epsilon\, $. The correponding completely isometric representation of 
$\, \mathfrak X \subset I ( \mathfrak X )\, $ thereby becomes a unital and $\epsilon$-unital $*$-representation with respect to superinvolution in $\,\mathcal B ( \widehat{\mathcal H} )\, $, hence completely $\epsilon$-positive. Consider its enveloping $\,\mathbb Z_2$-graded operator system $\,\widehat{\mathfrak X}\, $. To get the last statement of (i) recall that any 
$\mathbb Z_2$-graded operator system can be represented in standard form. Then any element is of the form 
$$ \begin{pmatrix} a & b \\ b & a \end{pmatrix} \> . $$ 
It is clear that 
$\,\Vert a \Vert = \inf_{b \in \widehat{\mathfrak X}_1} \Vert a + b \Vert\, $ and in the same manner 
$\, \Vert b \Vert = \inf_{a \in \widehat{\mathfrak X}_0} \Vert a + b \Vert\, $. So the natural projections 
$\, p_0\, $ and $\, p_1\, $ are complete quotient mappings. This proves (i). For the first statement in (ii) recall that by the Blecher-Ruan-Sinclair-Theorem 
$\, A\, $ has a unital completely isometric representation on Hilbert space, in particular $\, A\, $ being a unital operator space its enveloping operator system $\,\widehat{\mathfrak X}\, $ (naturally 
$\mathbb Z_2$-graded by composing the antilinear completely antiisometric involution on $\, A\, $ with the antilinear completely antiisometric adjoint operation) is well defined with injective envelope 
$\, I ( \widehat{\mathfrak X} , \bf 1 )\, $ a $\mathbb Z_2$-graded $C^*$-algebra. It is clear from the definition of multiplication on the injective envelope that the injection 
$\, A \subseteq \widehat{\mathfrak X} \subseteq I ( \widehat{\mathfrak X} , \bf 1 )\, $ is a $*$-homomorphism. To get our representation consider any unital graded $*$-representation of 
$\, I ( \widehat{\mathfrak X} , \bf 1 )\, $. Its restriction to $\, A\, $ will give the desired representation
\qed
\par\bigskip\noindent
The following result is well known (c.f. \cite{E-R}, Corollary 5.2.3) and included only for convenience of the reader since we will frequently use it in the sequel.
\par\bigskip\noindent
{\bf Proposition 1.}\quad Let $\, \Phi : A \rightarrow B\, $ be a bijective unital, linear complete isometry of two $C^*$-algebras $\, A\, $ and $\, B\, $. Then $\,\Phi\, $ is  a $*$-isomorphism\qed
\par\bigskip\noindent
The draw-back of the ordering by $\epsilon$-positive elements is that it is not intrinsic in the 
$*$-algebraic structure of a given super-$C^*$-algebra $\, A\, $. In practice the completed super operator system of $\, A\, $ will only be determined by a particular representation on graded Hilbert space. 
Here the $\epsilon$-order is of some value in distinguishing between different representations as demonstrated in the example preceding the representation theorem. We will therefore call two 
completely isometric $*$-representations of an incomplete super operator system 
{\it $\epsilon$-equivalent} iff the induced complete isometry of the represented concrete super operator systems extends to an $\epsilon $-unital complete isometry on the completed super operator systems. We will now sketch how to derive an intrinsic ordering which is independent of the representation. Thinking of how positive elements in a 
$C^*$-algebra are defined one may be tempted to define a matrix order on a given super-$C^*$-algebra $\, A\, $ by considering the closure $\, {\mathcal C}^*_{n , +} ( A )\, $ of the convex cone of positive linear combinations of elements of the form 
$\, x x^*\, $ with $\, x \in M_n ( A )\, $. In certain cases this will indeed result in a matrix order on the given 
super-$C^*$-algebra. This means that the resulting cones should linearly generate the sets of hermitian elements in $\, M_n ( A )\, $, such that $\, {\mathcal C}^*_{n , +} \cap - {\mathcal C}^*_{n , +} = 
\{ 0 \}\, $, and that for any hermitian element $\, x\, $ there exists a positive number $\, r \geq 0\, $ with 
$\, x + r {\bf 1}_n \in {\mathcal C}^*_{n , +}\, $. Also for a given matrix 
$\, \alpha \in M_{n , r} ( \mathbb C )\, $ one must have $\, \alpha\cdot {\mathcal C}^*_{r , +}\cdot \overline\alpha \subseteq {\mathcal C}^*_{n , +}\, $. A super $C^*$-algebra posessing this property 
will be called a {\it pure} super $C^*$-algebra. 
An example is the algebra $\, A_0\, $ generated by a single selfadjoint element 
$\, x\, $ with $\, x^2 = 0\, $ of the example above. The details are left to the reader. There are however also rather striking counterexamples. The simplest is given by the $\mathbb Z_2$-graded Clifford algebra $\, C_1\, $ on one odd generator. In this case the cones $\, {\mathcal C}^*_{n , +}\, $ make up all of the hermitian elements. The same holds for the $\mathbb Z_2$-graded $C^*$-algebra 
$\, \widehat M_2 ( B )\, $ with standard off-diagonal grading for a trivially graded $C^*$-algebra $\, B\, $.
If we are not to stick with the case of pure super-$C^*$-algebras, a property which in general is hard to check we must find some other way of introducing an intrinsic order. A key observation in this context is given by the proof of the following proposition.
\par\bigskip\noindent
{\bf Proposition 2.}\quad To any (unital) $\mathbb Z_2$-graded $C^*$-algebra $\,\widehat A\, $ one can associate a graded $C^*$-algebra $\, M_{\widehat A}\, $ which is a nontrivial continuous field of $C^*$-algebras over the unit circle with each fibre $*$-isomorphic to $\,\widehat A\, $ (the Moebius strip over 
$\,\widehat A\, $). 
\par\bigskip\noindent
{\it Proof.}\quad We will introduce a continuous family of new products on $\,\widehat A\, $, indexed by the points of the unit circle $\, S^1 = \{ \omega\in\mathbb C\,\vert\, \vert \omega\vert = 1 \}\, $ and show that the resulting continuous field $\, A\, $ when equipped with a suitable matrix norm becomes a unital $C^*$-algebra. 
Let elements $\, x = x_0 + x_1\, $ and 
$\, y = y_0 + y_1\, $ of $\, \widehat A\, $ be given where always $\, x_0\, $ denotes the even part of 
$\, x\, $ and $\, x_1\, $ the odd part. Define the ${*}_{\omega }$-product of $\, x\, $ and $\, y\, $ by the formula 
$$ x\, {*}_{\omega} y = ( x_0 y_0 + \omega x_1 y_1 ) + ( x_0 y_1 + x_1 y_0 ) \> . $$
Associativity of the product is readily checked. Since $\, M_n ( \widehat A )\, $ is also a 
$\mathbb Z_2$-graded $C^*$-algebra it is clear that  the $*$-multiplication extends to matrix algebras in a compatible way. Consider the algebra which is linearly isomorphic with 
$\, C ( S^1 )Ê\otimes \widehat A\, $ but with product at the fibre over the point $\,\omega\, $ twisted according to the above formula. It is clear that this local definition of the product coherently defines a global product on the tensor product $\, C ( S^1 ) \otimes \widehat A\, $ which is denoted $\, *\, $ to distinguish it from the ordinary (tensor) product. Let $\, x = x ( \omega )\, $ be an element of this algebra.  In order to obtain a $C^*$-algebra we will have to renormalize the fibres according to the following scheme. For each $\, \omega\, $ the fibre 
$\, A_{\omega }\, $ which is the image of our continuous field under the evaluation homomorphism corresponding to the point $\,\omega\, $ is isomorphic with $\,\widehat A\, $ by the map 
$$ {\iota }_{\omega }\> :\> A_{\omega } \>\longrightarrow \widehat A\quad , \quad 
x_0\> +\> x_1\>\mapsto\> x_0\> +\> \sqrt{\omega }\, x_1 $$
where the sign of $\,\sqrt{\omega }\, $ is deliberate (but fixed in order that $\, \sqrt{ \omega }\, $ is continuous in the parameter $\,\omega\, $). 
One may than use the map $\, {\iota }_{\omega }\, $ to renormalize the fibres $\, A_{\omega }\, $ such that each becomes isometric and $*$-isomorphic with $\,\widehat A\, $. The (fibrewise) involution and grading make this continuous field into a (graded) $C^*$-algebra. Note however that there is no global trivialization which identifies the continuous field with $\, C ( S^1 ) \otimes \widehat A\, $ 
\qed
\par\bigskip\noindent
In particular consider the twisted product at the point $\, \omega = -1\, $ together with the corresponding map 
$$ \iota\> :\> ( \widehat A\, ,\, *_{-1} )\>\longrightarrow\> ( \widehat A\, ,\, *_1 )\quad ,\quad 
x_0\, +\, x_1\>\mapsto\> x_0\, +\> i\, x_1 \> . $$
Let $\, A\, $ be an abstract unital super-$C^*$-algebra with hermitian subspace $\, A^h\, $. The 
enveloping $\,\mathbb Z_2$-graded operator system of $\, A\, $ with injective envelope 
equal to $\, I ( A )\, $ (which is a graded $C^*$-algebra) is then well defined, so is its 
enveloping (minimal $\mathbb Z_2$-graded) $C^*$-algebra $\,\widehat A\, $ inside this injective envelope.  Upon utilizing the map $\,\iota\, $ (viewed as a linear map of $\,M_n ( \widehat A )\, $ to itself) the hermitian elements of $\, M_n ( A )\, $ are identified with a unital subspace of 
the ordinary selfadjoint elements of $\, M_n ( \widehat A )\, $ which is the selfadjoint part of an operator system. Pulling back the order structure of the latter to $\, M_n ( A )\, $ defines a matrix order on $\, A\, $ with corresponding cones $\, {\mathcal C}^s_{n , +} \subseteq M_n ( A )^h\, $ of {\it superpositive elements}. If 
$\, \mathfrak X\, $ is any unital subspace of $\, A\, $ invariant under the superinvolution, it inherits a matrix order by taking the intersection of the cones of superpositive elements of $\, M_n ( A )\, $ with 
$\, M_n ( \mathfrak X )\, $. At first sight the definition of superpositive elements seems to be a simple trick, 
exchanging the subspaces of hermitian and ordinary selfadjoint elements in a $\mathbb Z_2$-graded frame by the linear map $\,\iota\, $. What is not so clear is how the norms of an element 
$\, x\, $ and its image $\, \iota ( x )\, $ are related (precisely). One easily derives that
$$  1 / 2\cdot \Vert x \Vert \leq \Vert \iota ( x ) \Vert \leq  2\cdot \Vert x \VertÊ\> , $$
which seems to be the best possible estimate in general. We wish to have an abstract notion of super operator system which is independent of any particular representation or embedding into a complete super operator system. One possibility is to follow the scheme of abstract operator systems by specifying a matrix order on the set of hermitian elements. In the case of an operator system the norm of an element is then uniquely determined by this order if one looks for a completely isometric {\it and} completely order isomorphic representation. The proposition below states that with super operator systems things are more complicated (the norm of an element is determined only within a global margin). Also, in practice,  the order is usually induced by a representation and is difficult to build in an abstract setting. We will therefore choose a different approach.
\par\bigskip\noindent
{\bf Definition 2.}\quad An (abstract) super operator system $\, ( \mathfrak X\, ,\, *\, ,\, {\bf 1} )\, $ is an abstract unital operator space $\, ( \mathfrak X\, ,\, {\bf 1} )\, $ (abstract operator space with specified unit element that admits a completely isometric unital representation on Hilbert space) equipped with a completely antiisometric antilinear involution $\, *\, $ fixing the unit. 
\par\smallskip\noindent
An (abstract) super operator space $\, ( \mathfrak V\, ,\, * )\, $ is an abstract operator space 
$\,\mathfrak V\, $ with a completely antiisometric antilinear involution $\, *\, $.
\par\smallskip\noindent
Given an operator space 
$\, \mathcal V\, $ with to different matrix norms $\,\{ \Vert\cdot {\Vert}^1_n {\}Ê}_{n\in\mathbb N}\, $ and 
$\,\{ \Vert\cdot {\Vert }^2_n {\} }_{n\in\mathbb N}\, $, these will be called {\it completely equivalent} iff there exists a global constants $\, C_1\, ,\, C_2 > 0\, $ with 
$$ C_1\cdot \Vert x {\Vert }^1_n \leq 
\Vert x {\Vert }^2_n \leq C_2\cdot \Vert x {\Vert }^1_n  $$
for all $\, n\, $. 
\par\bigskip\noindent
{\bf Theorem.}\quad (Representation theorem for super operator systems)
\par\noindent 
Let $\, \mathfrak X\, $ be an abstract super operator system. There exists a $\mathbb Z_2$-graded Hilbert space $\,\widehat{\mathcal H}\, $ and a completely isometric unital 
$*$-representation $\, \varphi : \mathfrak X \rightarrow \mathcal B ( \widehat{\mathcal H} )\, $. 
\par\smallskip\noindent
Let $\, \mathfrak V\, $ be an abstract super operator space. Then there exists a Hilbert space 
$\,\mathcal H\, $ and a completely isometric $*$-representation 
$\,\nu : \mathfrak V \rightarrow \mathcal B ( \mathcal H )\, $ of $\,\mathfrak V\, $ as a pre-operator system.
\par\bigskip\noindent
{\it Proof.}\quad Let $\,\psi : \mathfrak X \rightarrow \mathcal B ( \mathcal H )\, $ be a completely isometric unital representation of the unital operator space $\,\mathfrak X\, $. Define the representation  
$\  \overline{\psiÊ} : \mathfrak X \rightarrow \mathcal B ( \mathcal H )\, $ by composing $\,\psi\, $ with the 
$*$-involution of $\,\mathfrak X\, $ on the right and the adjoint operation on 
$\,\mathcal B ( \mathcal H )\, $ on the left. Then $\,\overline{\psi }\, $ is another completely isometric unital representation of $\,\mathfrak X\, $. Put $\,\widehat{\mathcal H } = \mathcal H \oplus \mathcal H\, $ and define a grading on the direct sum by exchanging the two copies of $\,\mathcal H\, $. It is easily verified that the direct sum representation $\, \psi \oplus \overline{\psi }\, $ has the required properties. Much in the same way one proves the existence of a completely isometric $*$-representation $\,\nu\, $ on a graded Hilbert space in case of a super operator space. Now given a completely isometric $*$-representation $\,\widehat\nu\, $ with respect to superinvolution on graded Hilbert space $\,\widehat{\mathcal H}\, $ let 
$\,\kappa\, $ be a square root of the grading operator $\,\epsilon\, $ as above and put 
$\, \nu\, =\, \kappa\cdot \widehat\nu\cdot\kappa\, $. Then $\,\nu\, $ defines a completely isometric $*$-representation of $\,\mathfrak V\, $ with respect to the ordinary (adjoint) involution on 
$\,\mathcal B ( \widehat{\mathcal H} )\, $ \qed 
\par\bigskip\noindent
One may now {\it define} the matrix order on $\, \mathfrak X\, $ by pulling back the cones of 
superpositive elements in the concrete super operator system $\,\varphi ( \mathfrak X )\, $. The Proposition below shows that this does not depend on the particular representation. Thus the structure of abstract super operator system uniquely determines the matrix order on $\,\mathfrak X\, $. On the other hand given the order there might exist different matrix norms compatible with the matrix order. Any two such matrix norms will then have to be completely equivalent from the following Proposition. 
\par\bigskip\noindent
{\bf Proposition 3.}\quad Let $\, \mathfrak X\, ,\, \mathfrak Y\, $ be super operator systems. If 
$\, \varphi : \mathfrak X \rightarrow \mathfrak Y\, $ is a unital and (completely) contractive $*$-linear map, then it is (completely) superpositive. A completely superpositive map 
$\,\psi : \mathfrak X \rightarrow \mathfrak Y\, $ is completely bounded with 
$\, \Vert \psi {\Vert }_{c b} \leq 8\cdot  \Vert \psi ( \bf 1 ) \Vert\, $. If $\psi\, $ is unital then 
$\, \Vert \psi {\Vert }_{c b} \leq 4\, $.
\par\bigskip\noindent
{\it Proof.}\quad Let $\,\varphi\, $ be a unital, contractive $*$-linear map. Then it extends to a 
graded contractive map on the enveloping graded operator systems 
$\,\widehat\varphi : \widehat{\mathfrak X} \rightarrow \widehat{\mathfrak Y}\, $. One has 
$\, \varphi ( x ) = \iota\circ \widehat\varphi\circ{\iota }^{-1} ( x )\, $ and since $\,\widehat\varphi\, $ is positive, $\,\varphi\, $ is superpositive. Let $\, \psi\, $ be completely superpositive. Then the map 
$\, {\iota }\circ\psi\circ {\iota }^{-1} : \iota ( \mathfrak X ) \rightarrow \iota ( \mathfrak Y )\, $ (viewed as subspaces of $\,\widehat{\mathfrak X}\, ,\, \widehat{\mathfrak Y}\, $) is completely positive, so that its completely bounded norm is given by $\, \Vert \iota\circ\psi\circ {\iota }^{-1} ( {\bf 1} ) \Vert \leq 
2\cdot \Vert \psi ( {\bf 1} ) \Vert\, $. Then the map $\,\psi\, $ must be completely bounded with completely bounded norm less or equal to $\, 8\cdot \Vert \psi ( {\bf 1} ) \Vert\, $. If $\,\psi\, $ is unital then also $\, \iota\circ\psi\circ {\iota }^{-1}\, $ whence the result\qed
\par\bigskip\noindent
Given the matrix order on an abstract super operator system one may consider the set of equivalence classes of all {\it completely superpositive unital} representations as opposed to completely contractive. Any such representation is completely bounded by $\, 4\, $. Let $\, \max \mathfrak X\, $ denote the image of $\,\mathfrak X\, $ under the direct sum of all such representations. Then $\, \max \mathfrak X\, $ is a super operator system with the same underlying vector space, unit and involution but equipped with a larger though completely equivalent matrix norm. It is clear from construction that $\,\max \mathfrak X\, $ has the additional property that any unital and completely superpositive map into another super operator system is completely contractive. A super operator system with this property will be called {\it maximal}. If one considers the {\it folium} of all super operator systems with underlying vector space $\,\mathfrak X\, $ and same matrix order (but different matrix norms) it contains a unique element which is an operator system (namely the image of $\, \mathfrak X\, $ in some completely isometric $*$-representation under the map $\,\iota\, $). This will be denoted by $\,{\mathfrak X}_0\, $. Given a super operator system 
$\,\mathfrak X\, $ one may consider its enveloping $\mathbb Z_2$-graded $C^*$-algebra in the direct sum of all (up to equivalence) unital completely contractive $*$-representations on graded Hilbert space. This object will be called the {\it maximal} enveloping graded $C^*$-algebra of $\,\mathfrak X\, $, denoted 
$\, C^*_{max} ( \mathfrak X )\, $. This definition slightly differs from the previously defined maximal enveloping graded $C^*$-algebra of a super-$C^*$-algebra $\, A\, $ which only takes into account (multiplicative) $*$-representations of the operator algebra. As before the construction is functorial in the following sense: given any completely contractive unital $*$-linear map of super operator systems 
$\, \mathfrak X \rightarrow \mathfrak Y\, $ there is a canonically associated graded 
$C^*$-homomorphism $\, C^*_{max} ( \mathfrak X ) \rightarrow C^*_{max} ( \mathfrak Y )\, $. 
If $\,\mathfrak V\, $ is a super operator space one defines similarly the {\it universal graded representation of 
$\,\mathfrak V\, $} by taking the direct sum representation over all (up to equivalence) completely contractive $*$-representations of $\,\mathfrak V\, $ on graded Hilbert space. 
Let $\,\widehat{\mathfrak V}\, $ denote the enveloping $\mathbb Z_2$-graded pre-operator system of 
$\, \mathfrak V\, $ in this representation. One defines $\, C^*_{max} ( \mathfrak V )\, $ to be the graded 
$C^*$-algebra generated by $\, \widehat{\mathfrak V}\, $. Again this definition slightly differs from the case of a super operator system, but if $\, \mathfrak V\, $ happens to be completely $*$-isometric with a super operator system $\, \mathfrak X\, $ then one has a natural graded surjective $C^*$-homomorphism $\, C^*_{max} ( \mathfrak V ) \twoheadrightarrow C^*_{max} ( \mathfrak X )\, $ and the construction is functorial in the usual sense. Since any super operator space can be represented as a pre-operator system on ungraded Hilbert space one may also consider the universal representation as a pre-operator system which is defined analogously. Let $\, C^*_{max , 0} ( \mathfrak V )\, $ denote the enveloping (ungraded) $C^*$-algebra of the image of $\,\mathfrak V\, $ in this representation. Also this construction is functorial.
For any element $\, x\, $ of a graded $C^*$-algebra $\,\widehat A\, $ define the graded absolute value of $\, x\, $ to be the unique superpositive square root $\, \vert x {\vert }_s\, $ of $\, x^**x\, $ making use of the $*$-multiplication (resp. the map $\,\iota\, $) as in the proof of Proposition 2. The graded absolute value is a special case of {\it graded continous function calculus} which can be applied to any element of 
$\,\widehat A\, $ which is normal with respect to $*$-multiplication. 
\par\bigskip\noindent
{\it Example.}\quad To see that super operator spaces arise quite naturally consider the dual operator space $\, X^*\, $ of an operator system $\, X\, $. The involution on $\, X\, $ determines an isometric involution on 
$\, X^* = \mathcal C\mathcal B ( X\, ,\, \mathbb C )\, $ by the formula 
$\, \overline f ( x ) = \overline{ f ( x^* )}\, $ (we use the notation $\,\overline f\, $ rather than $\, f^*\, $ to avoid mistaking it with the dual map $\, f^* \in \mathcal C\mathcal B ( {\mathbb C}^*\, ,\, X^* )\, $). 
If $\, \phi\, $ is an element of 
$\, M_n ( X^* ) = \mathcal C\mathcal B ( X\, ,\, M_n ( \mathbb C ) )\, $ then its completely bounded norm is given by the ordinary norm of the $n$-fold amplification 
$\,  {\phi }_n = {id}_{M_n ( \mathbb C )} \otimes \phi\, $ (compare \cite{Pa}, Prop. 8.11) from which it is easy to induce that the involution as above is in fact completely antiisometric.
Thus the dual $\, X^*\, $ is naturally a super operator space. 
\par\bigskip\bigskip\noindent 
{\bf 2. Strong matrix norms}
\par\bigskip\noindent 
Let $\, \mathfrak V\, ,\, \mathfrak W\, $ be super operator spaces. 
If $\, \phi : \mathfrak V \rightarrow \mathfrak W\, $ is a $*$-linear completely bounded map we will say that it is {\it hermitian $n$-contractive}, if the restriction of $\, {\phi }_n\, $ to the hermitian elements 
$\, M_n ( \mathfrak V )^h\, $ is contractive. One also defines the {\it hermitian $n$-norm} of a $*$-linear completely bounded map of super operator spaces by 
$$ \Vert \phi {\Vert }_n^h\> =\> \sup\, \{ \Vert {\phi }_n ( x ) \Vert\,\vert\, x \in M_n ( \mathfrak V )\, ,\, 
x = x^*\, ,\, \Vert x \Vert \leq 1 \} \> . $$
As with the ordinary norms it follows that one has 
$\, \Vert \phi {\Vert }_n^h = \Vert {\phi }^* {\Vert }_n^h\, $ for each $\, n\, $ where 
$\, {\phi }^* : {\mathfrak W}^* \rightarrow {\mathfrak V}^*\, $ is the dualized map. Finally, let 
$\, \mathcal V \subseteq \mathcal B ( \widehat{\mathcal H} )\, $ be a concrete operator space represented on a $\,\mathbb Z_2$-graded Hilbert space with grading operator $\,\epsilon\, $ (in case of an ungraded Hilbert space put $\, \epsilon = {\bf 1}\, $). For $\, x\in \mathcal V\, $ put 
$$ \Vert x {\Vert }^s\> =\> \sup\, \bigl{\{ } \vert \langle x \xi\, ,\, \epsilon \xi \rangle \vert\, \bigm\vert\, \xi \in \widehat{\mathcal H}\, ,\, 
\Vert \xi \Vert \leq 1 \bigr{\} } \> . $$
The resulting norm on $\, \mathcal V\, $ will be called the {\it strong norm} (with respect to the representation). One checks that if $\, \mathcal V\, $ is a super operator space and the representation is an isometric $*$-representation, then the involution is also isometric for the strong norm. Another nice property of the strong norm is that it is invariant under conjugation by unitaries of degree zero, i.e. if 
$\, U \in \mathcal B ( \widehat{\mathcal H} )_0\, $ is a unitary then 
$\, \Vert U\, x\, U^* {\Vert }^s = \Vert x {\Vert }^s\, $. More generally, one defines the strong $n$-norm on 
$\, M_n ( \mathcal V )\, $ to be its strong norm with respect to the natural identification 
$\, M_n ( \mathcal V ) \subseteq M_n ( \mathcal B ( \widehat{\mathcal H} ) ) \simeq \mathcal B ( {\widehat{\mathcal H} }^n )\, $. The sequence of strong norms $\, \{ \Vert\cdot {\Vert }_n^s \}\, $ satisfies axiom 
$$ \Vert x \oplus y {\Vert }_{n+m}\> =\> \max\, \{ \Vert x {\Vert }_n\, ,\, \Vert y {\Vert }_m \} \leqno{ \bf{(M 1)}}
$$
of a matrix norm, but instead of 
$$  \Vert \alpha\cdot x \cdot \beta {\Vert }_n\>\leq\> \Vert \alpha \Vert \, 
\Vert x {\Vert }_r\, \Vert \beta \Vert \> , \leqno{ \bf{(M 2)}} $$
for $\, \alpha\in M_{n , r} ( \mathbb C )\, ,\, \beta\in M_{r , n} ( \mathbb C )\, $ one only has the weaker property 
$$ \Vert \alpha\cdot x\cdot {\overline\alpha } {\Vert }_n^s\> \leq\> \Vert \alpha {\Vert }^2\cdot \Vert x {\Vert }_r^s \leqno{ {\bf (\Sigma M 2)} } $$ 
where $\, \alpha \in M_{n , r} ( \mathbb C )\, $ and $\, x\in M_r ( V )\, $. A sequence of norms satisfying 
$\,  {\bf (\Sigma M 1)}\, =\, {\bf (M 1)}\, $ and $\, {\bf (\Sigma M 2)}\, $ will be called a {\it strong matrix norm}. Given a Banach space 
$\, \mathcal V\, $ equipped with a strong matrix norm $\, \{ \Vert\cdot {\Vert }^s_n \}\, $ for the sequence of spaces 
$\,\{ M_n ( \mathcal V ) \}\, $ putting 
$$ \Vert x {\Vert }_n\> :=\> 2\, \Bigm\Vert \begin{pmatrix} 0 & x \\ 0 & 0 \end{pmatrix}Ê{\Bigm\Vert }^s_{2n} 
$$
defines an (ordinary) matrix norm on $\, \mathcal V\, $. Namely, suppose that 
$\, x = \alpha\cdot y\cdot\beta\, $ with $\, y\in M_r ( \mathcal V )\, ,\, \alpha\in M_{n , r} ( \mathbb C )\, ,\, 
\beta\in M_{r , n} ( \mathbb C )\, $ and without loss of generality that $\,\VertÊ\alpha \Vert \geq 
\Vert \beta \Vert\, $. Put $\, \widetilde\beta = \Vert \alpha \Vert\cdot \Vert \beta {\Vert }^{-1}\cdot \beta\, $ and $\,\widetilde y = \Vert \beta \Vert\cdot \Vert \alpha {\Vert }^{-1}\cdot y\, $. Then 
$$ \begin{pmatrix} 0 & x \\ 0 & 0 \end{pmatrix}\> =\> \begin{pmatrix} \alpha & 0 \\ 0 & \overline{\widetilde\beta} \end{pmatrix}\cdot \begin{pmatrix} 0 & \widetilde y \\ 0 & 0 \end{pmatrix}\cdot 
\begin{pmatrix} \overline\alpha & 0 \\ 0 & \widetilde\beta \end{pmatrix} $$ 
so that 
$$ \Vert x {\Vert }_n\> \leq\> \Vert \alpha \Vert\cdot \Vert \widetilde y {\Vert }_r\cdot \Vert \widetilde\beta \Vert\> =\> \Vert \alpha \Vert\cdot \Vert y {\Vert }_r\cdot \Vert \beta \Vert\> . $$
Property $\, {\bf ( M1)}\, $ is even more trivial to check. Note that in any case the following relation holds 
$$ \Vert x {\Vert }^s_n\> =\> \Bigm\Vert \begin{pmatrix} x & 0 \\ 0 & -x \end{pmatrix} {\Bigm\Vert }^s_{2n}\> =\> \Bigm\Vert \begin{pmatrix} 0 & x \\ x & 0 \end{pmatrix} {\Bigm\Vert }^s_{2n} \> \leq\> 
2\,\Bigm\Vert \begin{pmatrix} 0 & x \\ 0 & 0 \end{pmatrix} {\Bigm\Vert }^s_{2n} \> =\> 
\Vert x {\Vert }_n\> . $$
It is also clear that in case of a concete operator space the matrix norm thus obtained from the strong matrix norm coincides with the matrix norm of the representation whereas there may exist representations of the same underlying matrix norm giving rise to different strong norms. Thus a strong matrix norm imposes a {\it finer} structure on an operator space than an ordinary matrix norm. This leaves unanswered the question which abstract strong matrix norms can actually be represented. In case of a super operator space $\,\mathfrak V\, $ and a $*$-representation on (graded) Hilbert space one will certainly need at least one more axiom besides demanding that the involution is strongly completely antiisometric. In particular for each hermitian 
element $\, x  \in M_n ( \mathfrak V )^h\, $ one must have 
$$ \Vert x {\Vert }_n^s\> =\> 2\, \Bigm\Vert \begin{pmatrix} 0 & x \\ 0 & 0 \end{pmatrix} {\Bigm\Vert }_{2n}^s\quad ,\quad x \in M_n ( \mathfrak V )^h 
\leqno{ {\bf (h\Sigma M*)}} $$ 
since for a selfadjoint (diagonalizable) element the strong norm coincides with the usual norm.  However since this property is valid only for hermitian elements it does not seem sufficient to guarantee the existence of a strongly completely isometric $*$-representation. 
 A linear map 
$\,\phi : \mathcal V \rightarrow \mathcal W\, $ of 
concrete operator spaces will be called {\it (completely) strongly contractive} if it is (completely) contractive with respect to the strong norms. It will be called {\it really strongly contractive} if for each unit vector $\, \eta\, $ in the representation space of $\, \mathcal W\, $ and each $\, x\in \mathcal V\, $ there exists a unit vector 
$\, \xi\, $ in the representation space of 
$\, \mathcal V\, $ such that putting $\, \alpha = \langle x\,\xi\, ,\, \epsilon\,\xi\rangle\, $, 
$\, \beta = \langle \phi ( x )\,\eta\, ,\, \epsilon\,\eta \rangle\, $ one has the relations
$$ \vert \alpha \vert\> \geq\> \vert \beta \vert\quad ,\quad 
 \vert Re ( \alpha )\vert \>\geq\>\vert Re ( \beta ) \vert \> . $$
The first condition just states that $\,\phi\, $ is strongly contractive. As for the second it looks a bit awkward and difficult to check in a general situation but we will see some instances in finite dimensional matrix algebras where it can be sucessfully applied. If with notation as above one also has that 
the sign of $\, Re ( \alpha )\, $ is the same as the sign of $\, Re ( \beta )\, $ the map is called 
{\it signed really strongly contractive}.
\par\bigskip\noindent
{\bf Strong Lemma. }\quad Let $\, \phi : \mathfrak V \rightarrow \mathfrak W\, $ be a $*$-linear map of super operator spaces and $\, \mathcal T \subseteq \mathfrak V\, $ a closed subspace with 
$\, \mathcal T + {\mathcal T}^* = \mathfrak V\, $. If there exist isometric $*$-representations of 
$\,\mathfrak V\, $ and 
$\, \mathfrak W\, $ such that $\, \phi\, $ restricted to 
$\, \mathcal T\, $ is really strongly contractive, then $\,\phi\, $ is strongly contractive (and hence hermitian contractive).
\par\bigskip\noindent
{\it Proof.}\quad  If $\,\mathfrak V \subseteq \mathcal B ( \widehat{\mathcal H} )\, $ consider  the strongly isometric cross-diagonal embedding 
$\, \lambda\, :\, \mathfrak V \rightarrow M_2 ( \mathcal B ( \widehat{\mathcal H} ) )\, $ by the assignment
$$ x\> \mapsto\> \begin{pmatrix} 0 &  x  \\  x  & 0 \end{pmatrix} \> . $$
That it is strongly isometric follows from the fact that the image of $\, x\, $ can be rotated to 
$$ \begin{pmatrix}  x & 0 \\ 0 & - x  \end{pmatrix}\> .  $$
Now consider the operator subsystem $\, \widetilde{\mathfrak V}\, $ of 
$\, M_2 ( \mathcal B ( \widehat{\mathcal H} ) )\, $ obtained by adjoining the grading operator  (still denoted $\,\epsilon\, $ for simplicity) to 
$\, \lambda ( \mathfrak V )\, $. 
The same procedure applies to $\, \mathfrak W\, $. Also put 
$\, \widetilde{\mathcal T}\, =\, \lambda  ( \mathcal T ) + \mathbb C\cdot \epsilon\, $. 
We claim that the natural $\epsilon$-unital extension $\,\widetilde{\phi }\, $  of $\, \phi\, $  is strongly contractive when restricted to $\,\widetilde{\mathcal T}\, $. 
If $\, c\in\mathbb C\, $ is any complex number, then 
$$ \Bigm\Vert \begin{pmatrix} c\cdot\epsilon & \phi ( x )  \\ 
\phi ( x ) & c\cdot\epsilon\end{pmatrix} 
{\Bigm\Vert }^s\> =\>   \max \bigl{\{ } \Vert \phi ( x ) + c\cdot\epsilon {\Vert }^s \, ,\, 
\Vert \phi ( x ) - c\cdot\epsilon {\Vert }^s  \bigr{\} } \> . $$
Let $\, c = {\omega }_c\cdot \gamma\, $ with $\,\gamma \geq 0\, $ and $\, \vert {\omega }_c \vert = 1\, $.
Dividing by $\, {\omega }_c\, $ will not change anything about the strong norms so we may assume that 
$\, c = \gamma\, $ is a positive number. Let $\, \eta\, $ be a unit vector in the representation space of 
$\, \widetilde{\mathfrak W}\, $ such that the maximum of the strong norms of 
$\, \phi ( x ) \pm \gamma\, \epsilon\  $ is attained at $\,\eta\, $. Then there exists a unit vector $\, \xi\, $ in 
$\,\widehat{\mathcal H}\, $ such that 
$$\vert  \langle  x \, \xi\, ,\, \epsilon\, \xi \rangle \vert\> \geq\> 
\vert \langle \phi ( x )\, \eta\, ,\, \epsilon\,\eta \rangle \vert \>  $$
and 
$$ \vert Re ( \langle x\,\xi\, ,\, \epsilon\,\xi \rangle )Ê\vert\>\geq\> 
\vert Re ( \langle \phi ( x )\,\eta\, ,\, \epsilon\,\eta \rangle\vert \> . $$
Now if $\,\alpha\, $ and $\, \beta\, $ are complex numbers with 
$\,\vert \alpha \vert \geq \vert \beta \vert\, $ and $\, Re ( \alpha ) \geq Im ( \alpha )\, $ then also 
$$ \max\, \{ \vert \alpha + \gamma \vert\, ,\, \vert \alpha - \gamma \vert \}\> \geq\> 
\max\, \{ \vert \beta + \gamma \vert\, ,\, \vert \beta - \gamma \vert \} $$
for any positive number $\, \gamma\, $. This proves the claim. 
Let then $\, \eta\, $ be an arbitrary unit vector in the representation space of 
$\, \widetilde{\mathfrak W} \, $ and consider the functional 
$$ \rho ( x ) \> =\> \langle \lambda ( \phi ( x ) )\,\eta\, ,\, \epsilon\, \eta \rangle $$
on $\,\widetilde{\mathcal T}\, $ which is contractive for the strong norm. 
Then from the Hahn-Banach-Theorem it extends to an $\epsilon$-unital and contractive functional 
$\,\widetilde\rho\, $ on $\, \widetilde{\mathfrak V}\, $ with respect to the strong norm which means that the functional is hermitian contractive with respect to the ordinary norm. It is also contractive on the larger real subspace generated by the hermitian elements of 
$\, \widetilde{\mathfrak V}\, $ and complex multiples of $\, \epsilon\, $. To see this one utilizes the strong isometry 
$$ x\> \mapsto\> \kappa\cdot x\cdot\kappa\quad ,\quad \epsilon\>\mapsto\> {\bf 1} $$
with $\,\kappa\, $ a square root of the grading operator on $\,\mathcal B ( \widetilde{\mathcal H} )\, $ to reduce to the situation of linear combinations of selfadjoint elements in a pre-operator system and complex multiples of the identity. Any such element is diagonalizable so that
the strong norm coincides with its ordinary norm.
Now apply Lemma 2 above (resp. Lemma A.4.2 of \cite{E-R}) to conclude
that $\,\widetilde\rho\, $ is selfadjoint and hermitian contractive, by which one deduces that 
$$ \widetilde{ \rho } ( x )\> =\> \langle \widetilde{\phi } ( x ) \eta\, ,\, \eta \rangle $$
holds for all $\, x\in \widetilde{\mathfrak V}\, $. This means upon taking the supremum over all unit vectors 
$\,\eta\, $ that $\,\phi\, $ is hermitian contractive as stated\qed
\par\bigskip\noindent
{\it Example.}\quad Consider a onedimensional projection in $\, M_m ( \mathbb C )\, $, say with zero entries in all matrix entries except a $1$ in the lower right corner 
$$ p\> =\> \begin{pmatrix}   0\>\cdots\> 0 & 0 \\ \vdots\ddots\vdots & \vdots \\ 0\>\cdots\> 0 & 0 \\   
0\>\cdots\> 0 & 1 \end{pmatrix}\> . $$
and the projection $\,\varphi\, $ with kernel generated by $\, p\, $ and range the complementary subspace $\, \mathfrak Q\, $ of codimension one whose elements have zero entry in the lower right corner
$$ \mathfrak Q\> =\> \begin{pmatrix}   *\>\cdots\> * & * \\ \vdots\ddots\vdots & \vdots \\ *\>\cdots\> * & * \\   
*\>\cdots\> * & 0 \end{pmatrix}\> . $$
Let $\, \mathcal T\, $ be the subspace of lower triangular matrices in $\, M_m ( \mathbb C )\, $. Then the condition of the Sublemma above is easily seen to hold in this situation, so that $\,\varphi\, $ is hermitian contractive for any $\, m\, $.
As a second example consider the onedimensional subspace of  $\, M_m ( \mathcal C )\, $ spanned by the element
$$ e\> =\> \begin{pmatrix}   0\>\cdots\> 0 & 0 & 0 \\ \vdots\ddots\vdots & \vdots & \vdots \\ 
0\>\cdots\> 0 & 0 & 0 \\ 0\>\cdots\> 0 & 1 & 0 \\   
0\>\cdots\> 0 & 0 & -1\>\> \end{pmatrix}\>  $$
and the projection $\,\varphi\, $ with kernel generated by $\, e\, $ and range the subspace of codimension one whose entries look like
$$ \mathfrak F\> =\> \begin{pmatrix}   *\>\cdots\> * & * & * \\ \vdots\ddots\vdots & \vdots & \vdots \\ 
*\>\cdots\> * & * & * \\ *\>\cdots\> * & \alpha & \beta \\   
* \>\cdots\> * & \gamma & \alpha\>\> \end{pmatrix}\>  . $$
Again one checks the condition of the Sublemma for the lower triangular matrices with respect to this projection and concludes that it is hermitian contractive.
Similar examples of this kind can be found.
\par\smallskip\noindent
Besides the strong matrix norm as above there is another strong norm for the $C^*$-algebra of bounded operators on graded Hilbert space which is of interest and will be called the {\it strong $\sigma$-norm}. It is defined as follows. Let $\, x = x_0 + x_1\, $ be an operator on a graded Hilbert space 
$\,\widehat{\mathcal H}\, $ with grading operator $\,\epsilon\, $ such that $\, x_0\, $ is homogenous of degree zero and $\, x_1\, $ is of degree one. Define 
$$ \Vert x {\Vert }^{\sigma }\> :=\> \sup\, \{ \vert \langle x_0 \xi\, ,\, \xi \rangle\, +\, \langle x_1 \xi\, ,\, \epsilon \xi \rangle \vert\> \bigm\vert\> \xi \in \widehat{\mathcal H}\, ,\, \Vert \xi \Vert \leq 1  \} \> . $$
Similarly one defines the strong $\sigma $-$n$-norm to be the $\sigma $-strong norm of 
$\, \mathcal B ( {\widehat{\mathcal H}}^n )\simeq M_n ( \mathcal B ( \widehat{\mathcal H} ) )\, $.
Then it is clear how to define the strong matrix $\sigma $-norm for a concete (super) operator space 
$\,\mathcal V\, $ represented on $\,\widehat{\mathcal H}\, $ and one checks that the involution is completely antiisometric also with respect to the strong $\sigma $-norms. The strong matrix $\sigma $-norm still satisfies axioms $\, {\bf (\Sigma M 1)}\, $ and $\, {\bf (\Sigma M 2)}\, $, but in this case the matrix norm given by 
$$ \Vert x {\Vert }_{\sigma , n}\> :=\> 2 \Bigm\Vert \begin{pmatrix} 0 & x \\ 0 & 0 \end{pmatrix} {\Bigm\Vert }_{2n}^{\sigma } $$
need not coincide with the ordinary operator matrix norm of $\, x\in M_n ( \mathcal V )\, $. We may call it the {\it matrix $\sigma $-norm}. While the strong norm is connected with the concept of $\epsilon$-positivity, the strong $\sigma $-norm is related to superpositivity. Namely, one checks that an element 
$\, x\in \mathcal B ( \widehat{\mathcal H} )\, $ is 
$\epsilon $-positive if and only if 
$$ \langle x \xi\, ,\, \epsilon \xi \rangle \geq 0 $$
for all $\, \xi\in \widehat{\mathcal H}\, $. Similarly, $\, x = x_0 + x_1\, $ is superpositive if and only if 
$$ \langle x_0 \xi\, ,\, \xi \rangle\, +\, \langle x_1 \xi\, ,\, \epsilon \xi \rangle \geq 0 $$
for all $\, \xi\, $. To see this define the isometry $\,\kappa : \widehat{\mathcal H} \rightarrow \widehat{\mathcal H}\, $ by $\, \kappa ( {\xi }_0 \oplus {\xi }_1 ) = {\xi}_0 \oplus i {\xi }_1\, $. Now 
$\, x = x_0 + x_1\, $ is superpositive iff $\, \iota ( x ) = x_0 + i x_1\, $ is positive in the ordinary sense and one has the identity 
$$ \langle \iota ( x ) \xi\, ,\, \xi \rangle\> =\> 
\langle x_0\, \kappa \xi\, ,\, \kappa \xi  \rangle\, +\, \langle x_1\, \kappa \xi \, ,\, \epsilon\, \kappa \xi \rangle\> . $$
It is interesting to note that in case of a super operator system the matrix $\sigma $-norm, which is the ordinary matrix norm of $\,\iota ( \mathfrak X )\, $, is determined by its operator matrix norm (in particular it doesn't depend on the representation). One checks that the proof of the Strong Lemma applies just the same if instead of ordinary (strong) matrix norms (strong) matrix $\sigma $-norms are considered, so one only needs to replace the terms "really strongly contractive" and "hermitian contractive" of the Lemma by the corresponding terms "really strongly $\sigma $-contractive" and
"hermitian $\sigma $-contractive" respectively. 
\par\smallskip\noindent
Suppose given a super operator space $\, \mathfrak V\, $  equipped with a strong matrix norm 
$\,\{ \Vert \cdot {\Vert }^s_n {\}}_{n\geq 1}\, $ satisfying 
$\, {\bf (\Sigma M1)}\, $ and $\, {\bf (\Sigma M2)}\, $ and such that the involution is strongly completely antiisometric. Such a super operator space will be called an 
{\it abstract h-strong super operator space} if in addition it satisfies axiom $\, {\bf (h\Sigma M*)}\, $. One may define two other abstract strong matrix norms $\,\{ r_n\, ,\, R_n {\} }_{n\geq 1}\, $ from the given one by considering the maximum of the real matrix seminorms 
$$ p^{\omega }_n ( x )\> =\> {1\over 2}\, \Vert \omega x + \overline{\omega } x^* {\Vert }^s_n\quad ,\quad 
\omega\in\mathbb C\> ,\> \vertÊ\omega \vert = 1   $$
i.e $\, r_n ( x )\, =\, \sup_{\omega } \{ p^{\omega }_n ( x )\, \}\, $ and the maximum of the real matrix seminorms 
$$ P^{\omega }_n ( x )\> =\> {1\over 2}\,\Vert \omega x + {\overline\omega } x^* {\Vert }_n\quad ,\quad \omega\in\mathbb C ,\> \vert\omega\vert = 1 $$
i.e. $\, R_n ( x ) = \sup_{\omega } \{ P_n^{\omega } ( x ) \}\, $. 
If the original strong matrix norm is $*$-representable one gets 
$\, \Vert \omega x + {\overline\omega } x^* {\Vert }_n = \Vert \omega x + {\overline\omega } x^* 
{\Vert }^s_n\, $, and hence 
$$ r_n ( x )\quad =\quad R_n ( x )\quad \leq\quad \Vert x {\Vert }^s_n \> . $$
The reverse inequality can also be derived in this case. For given $\, x\, $ choose a complex number 
$\,\omega\, $ of modulus one such that $\, r_n ( x ) = 1/2\cdot\Vert \omega x + {\overline\omega } x^* 
{\Vert }^s_n\, $. Since the strong matrix norm is representable there exists a $*$-linear strongly contractive functional $\, \varphi  : M_n ( \mathfrak V ) \rightarrow \mathbb C\, $ with 
$\, \vert \varphi ( x ) \vert = \Vert x {\Vert }^s_n\, $
since in the induced representation of $\, M_n ( \mathfrak V )\, $ the strong norm of $\, x\, $ is given as the supremum over the absolute values of vector states 
$$ {\varphi }_{\xi } ( x )\> =\> \langle\, x \xi\, ,\, \xi\,\rangle\> .  $$
Assume that the inequality above is strict, so there exists a real number $\, 1/2 < r \leq 1\, $ with 
$\, \Vert x {\Vert }^s_n = r \Vert \omega x + {\overline\omega } x^* {\Vert }^s_n\, $. Identifying the range of 
$\,\varphi\, $ with the onedimensional subspace of $\, M_n ( \mathfrak V )\, $ generated by 
the element $\, \omega x + {\overline\omega } x^*\, $ assume that 
$\, \varphi ( x ) = \gamma r ( \omega x + {\overline\omega } x^* )\, $ with $\,\gamma\, $ a complex number of modulus one. Then 
$$ \varphi ( 1/2 ( {\overline\gamma } x + \gamma x^* ))\quad =\quad 
r ( \omega x + {\overline\omega } x^* ) $$
contradicting the assumption that $\,\varphi\, $ is strongly contractive. Since for any super operator space its ordinary matrix norm is $*$-representable from the representation theorem above this implies that a strong matrix norm on a super operator space is representable if and only if it satisfies the axiom 
$$ r_n ( x )\> =\> R_n ( x )\> =\> \Vert x {\Vert }^s_n \> . \leqno{ {\bf (\Sigma M*)}} $$
\par\medskip\noindent
A super operator space $\,\mathfrak V\, $ equipped with a strong matrix norm satisfying 
$\, {\bf (\Sigma M1)\, ,\, (\Sigma M2)}\, $ and $\, {\bf (\Sigma M*)}\, $ will be called an {\it abstract strong super operator space}. Note that the axiom also entails that the involution is completely antiisometric. 
The following theorem is a Corollary of the preceding discussion
\par\bigskip\noindent
{\bf Theorem.}$\> $ (Representation theorem for strong super operator spaces)
Let $\, \mathfrak V\, $ be an abstract strong super operator space. Then there exists an (ungraded) Hilbert space $\, \mathcal H\, $ and a strongly completely isometric $*$-representation 
$$ \mu : \mathfrak V\> \longrightarrow\> \mathcal B ( \mathcal H ) \> \qed $$
\par\bigskip\noindent
{\bf Theorem.}$\> $ (Representation of h-strong super operator spaces)
Let $\, \mathfrak V\, $ be an abstract h-strong super operator space. Then there exists for each $\, n\, $ an (ungraded) Hilbert space $\, {\mathcal H}_n\, $ and a strongly contractive and (strongly) hermitian isometric $*$-representation
$$ \nu : M_n ( \mathfrak V )\> \longrightarrow\> \mathcal B ( {\mathcal H}_n ) \> . $$
\par\bigskip\noindent
{\it Proof.}\quad  For the given 
h-strong super operator space replace the strong matrix norm by the abstract strong matrix norm 
$\, r_n ( x )\, =\, R_n ( x )\, $. 
It is easy to check that $\, r_n ( x )\, \leq\, \Vert x {\Vert }^s_n\, $ and that the induced matrix norms are the same, since $\, r_n ( x ) = R_n ( x )\, $ and
$$ \Vert x {\Vert }_{R_n}\> =\> 2\, R_{2n} \bigl( \begin{pmatrix} 0 & x \\ 0 & 0 \end{pmatrix} \bigr)\> =\> 
\sup_{\omega }\, \bigl\{ \Bigm\Vert \begin{pmatrix} 0 & \omega x \\ {\overline\omega } x^* & 0 \end{pmatrix} {\Bigm\Vert }_{2n} \bigr\}\> =\> \Vert x {\Vert }_n \> . $$
Since $\, \{ r_n \}\, $ is representable from the previous theorem this then also implies that 
$\, r_n ( x ) = \Vert x {\Vert }^s_n = \Vert x {\Vert }_n\, $ in case that $\, x = x^*\, $ is hermitian, since 
$$ \Vert x {\Vert }_n = \Vert x {\Vert }_{r_n} = r_n ( x ) \leq \Vert x {\Vert }^s_n $$
and the reverse inequality is trivial.
Suppose given $\, x\in M_n ( \mathfrak V )\, $ that the maximum of the set 
$\, \{ p^{\omega }_n ( x ) \}\, $ is taken at $\, \omega = 1\, $. From the Hahn-Banach Theorem one chooses a contractive projection $\,\varphi \, $ of $\, M_n ( \mathfrak V )\, $ onto the onedimensional subspace generated by $\, x + x^*\, $. Since the subspace is selfadjoint one can assume the projection to be $*$-linear, replacing it with $\, 1/2 ( \varphi + \overline{\varphi } )\, $ if necessary. Let $\, \varphi ( x ) = \alpha ( x + x^* )\, $ with $\,\alpha\in\mathbb C\, $. Since $\,\varphi\, $ is hermitian one gets 
$\, \alpha + \overline{\alpha } = 1\, $. If $\,\alpha\, $ is real then it must be $\, 1 / 2\, $ and one puts 
$\, f^1_x ( y ) = \varphi ( y )\, $. If on the other hand $\,\alpha\, $ is not a real number it follows that the maximum of the set $\, \{ p^{\omega }_n ( x ) \}\, $ cannot be taken at $\,\omega = 1\, $. 
One similarly constructs a $*$-linear function $\, f^{\omega }_x\, $ for each $\,\omega\, $ such that 
$\, \vert f^{\omega }_x ( y )\vert \leq r_n ( y )\, $ and $\, \vert f^{\omega }_x ( x )\vert  = p^{\omega }_n ( x )\, $ provided that the function $\, \sigma \mapsto p^{\sigma }_n ( x )\, $ attaines its maximum at the point 
$\, \omega\, $. Now take the direct sum over all onedimensional $*$-representations 
$\,\{ f^{\omega }_xÊ{\}Ê}_{x\in M_n ( \mathfrak V )} \, $ to get the result\qed
\par\bigskip\noindent
\par\bigskip\noindent
{\bf 3. Operator tensor products}
\par\bigskip\noindent
We now turn to the subject of tensor products of operator spaces 
resp. super operator systems. Since the operator spaces we are considering carry some extra structure we would like this structure to be reflected in the tensor product. For example, iff $\, \mathcal V\, $ and 
$\, \mathcal W\, $ are unital operator spaces it is natural to require that the tensor product 
$\, \mathcal V {\otimes }_{\gammaÊ} \mathcal W\, $ obtained as the completion of the algebraic tensor product by a certain cross matrix norm $\,\gamma\, $ also be unital with unit $\, 1_{\mathcal V} \otimes 1_{\mathcal W}\, $. We make the following convention: given an abstract operator space $\, \mathcal V\, $ together with a specified element $\, *_{\mathcal V} \in V\, $ of norm one 
(= {\it pointed operator space}), one defines the abstract unital operator space or {\it unification of}
$\, ( \mathcal V\, ,\, *_{\mathcal V} )\, $ to be completely isometric with the image 
$\, \overline{\mathcal V}\, $ of $\, \mathcal V\, $ under the direct sum of all completely contractive unital representations of $\, ( \mathcal V\, ,\, *_{\mathcal V} )\, $ on some Hilbert space up to equivalence. This space may or may not be completely isometric with $\, \mathcal V\, $. Since there always exists at least one completely contractive unital map into $\,\mathbb C\, $ by the Hahn-Banach theorem one sees that $\,\overline{\mathcal V}\, $ is not empty and is a unital operator space. 
If $\, \mathfrak V\, $ happens to be a super operator space with specified element 
$\, *_{\mathfrak V}\, $ of norm one then the pair $\, ( \mathfrak V\, ,\, *_{\mathfrak V} )\, $ is said to be a 
{\it pointed super operator space} if the involution fixes the element $\, *_{\mathfrak V}\, $. It is easy to see (compare with the proof of Lemma 4 below) that the unification of a pointed super operator space is a super operator system.
Given two unital operator spaces $\, ( \mathcal V\, ,\, 1_{\mathcal V} )\, $ and 
$\, ( \mathcal W , 1_{\mathcal W} )\, $ define the {\it unital projective tensor product} 
$\, \mathcal V \widehat{\otimes }_u \mathcal W\, $ to be the abstract unital operator space which is the unification of 
$\, ( \mathcal V \widehat\otimes \mathcal W \, ,\, 1_{\mathcal V} \otimes 1_{\mathcal W} )\, $ where 
$\, \mathcal V \widehat\otimes \mathcal W\, $ denotes the usual projective operator tensor product. The  
Haagerup tensor product $\, \mathcal V {\otimes}_h \mathcal W\, $ and the injective tensor product 
$\, \mathcal V \Check\otimes  \mathcal W\, $ of two unital operator spaces are known to be unital in any case but in case of the projective tensor product the author doesn't know if this is true. In any case it is easy to see that the unital projective tensor product is the maximal unital tensor product and surjects onto any other such tensor product in a canonical way. Let $\, ( \mathcal V\, ,\, *_{\mathcal V} )\, ,\, 
( \mathcal W\, ,\, *_{\mathcal W} )\, $ be pointed operator spaces. Then a completely bounded map 
$\, \varphi : \mathcal V \rightarrow \mathcal W\, $ is said to be {\it pointed} iff 
$\, \varphi ( *_{\mathcal V} ) = *_{\mathcal W}\, $. A unital operator space $\, \mathcal V\, $ is naturally a pointed operator space with $\, *_{\mathcal V} =  1_{\mathcal V} \, $. The mapping space of completely bounded linear maps of pointed operator spaces contains the specified subset of  completely contractive pointed maps denoted 
$\,\mathcal C\mathcal C\mathcal P ( \mathcal V\, ,\, \mathcal W )\, $. Let $\, \mathcal V\, ,\, \mathcal W\, $ and $\, \mathcal X\, $ be unital operator spaces. 
$\, \mathcal C\mathcal C\mathcal P\, ( \mathcal V \times \mathcal W \, ,\, \mathcal X )\, $ is the set of  pointed completely contractive {\it bilinear} maps $\,\varphi\, $ from $\, \mathcal V \times \mathcal W\, $ (with specified element $\, 1_{\mathcal V} \times 1_{\mathcal W}\, $) into 
$\, \mathcal X\, $, i.e. the maps are required to send 
$\, 1_{\mathcal V} \times 1_{\mathcal W}\, $ to $\, 1_{\mathcal X} \, $, and the c.b. norm is defined as usual in the mapping space of all bilinear maps by 
$\, \Vert \varphi {\Vert }_{cb} = \sup \{ \Vert {\varphi }_{p\times q; r\times s} \Vert\,\vert\, p , q , r , s \in \mathbb N \}\, =\, \sup \{ \Vert {\varphi }_{p \times p; p \times p}\,\vert\, p \in \mathbb N \}$,
$$ {\varphi }_{p , q; r , s} : M_{p , q} ( \mathcal V ) \times M_{r , s} ( \mathcal W ) \rightarrow 
M_{p \times r ; q\times s} ( \mathcal V \times \mathcal W )\> , $$
$$ {\varphi }_{p , q; r , s} ( [ v_{i , j} ]\, ,\, [ w_{k , l} ] ) = [ \varphi ( v_{i , j} , w_{k , l} ) ] $$
and the matrix norm is obtained using the identification 
$$ M_n ( \mathcal C\mathcal C\mathcal P ( \mathcal V \times \mathcal W\, ,\, \mathcal X ) ) \subseteq 
M_n ( \mathcal C\mathcal B\, ( \mathcal V \times \mathcal W\, ,\, \mathcal X ) ) = 
\mathcal C\mathcal B\, ( V \times W\, ,\, M_n ( X ) ) $$ 
(compare [Effros-Ruan], chap. 7).
\par\bigskip\noindent
{\bf Lemma 3.}\quad If $\, \mathcal V\, ,\, \mathcal W\, $ and $\, \mathcal X\, $ are unital operator spaces there are natural completely isometric identifications
$$ \mathcal C\mathcal C\mathcal P\, ( \mathcal V \widehat{\otimes }_u \mathcal W\, ,\, \mathcal X ) \simeq 
\mathcal C\mathcal C\mathcal P\, ( \mathcal V \times \mathcal W\, ,\, \mathcal X ) \simeq 
\mathcal C\mathcal C\mathcal P\, ( \mathcal V\, ,\, \mathcal C\mathcal B\, ( \mathcal W\, ,\, \mathcal X ) )\> . $$
\par\bigskip\noindent
{\it Proof}\quad Note that $\, \mathcal C\mathcal C\mathcal P\, ( \mathcal V \widehat\otimes \mathcal W\, ,\, \mathcal X )\,\simeq\, 
\mathcal C\mathcal C\mathcal P\, ( \mathcal V \widehat{\otimes }_u\, \mathcal W\, ,\, \mathcal X )\, $ 
with respect to the pointed projective tensor product $\, (  \mathcal V \widehat\otimes \mathcal W\, ,\, 
1_{\mathcal V} \otimes 1_{\mathcal W} )\, $ and compare with the proof of \cite{E-R}, Prop. 7.1.2  \qed
\par\bigskip\noindent 
The completely bounded maps of the unital projective tensor product are not that easily described, but the Lemma is sufficient for our purposes. Assume now that $\, \mathcal V\, $ is a unital operator space represented on some Hilbert space and let $\, \overline{\mathcal V}\, $ be the adjoint operator space. The mapping 
$\, v \mapsto \overline v\, $ is completely antiisometric in the sense of Definition 1. 
Note that in the case  where $\, \mathcal V\, $ is represented on a $\mathbb Z_2$-graded Hilbert space its ordinary adjoint is of course completely isometric with its superadjoint so there is nothing to be gained from the more general setting, but in turn one finds that the first case of the Lemma below also applies when $\, \mathcal V\, $ is a super operator space whence $\, \mathcal V\, $ is completely isometric with 
$\, \overline{\mathcal V}\, $.
Also let $\,\mathfrak X\, ,\, \mathfrak Y\, $ be two super operator systems with corresponding units and superinvolutions. 
\par\bigskip\noindent
{\bf Lemma 4.}\, The unital projective tensor product 
$\, \mathcal V \widehat{\otimes }_u \overline{\mathcal V}\, $, the Haagerup tensor product 
$\, \mathcal V {\otimesÊ}_h \overline{\mathcal V}\, $, and the injective tensor product 
$\, \mathcal V \check\otimes\overline{\mathcal V}\, $   of a unital operator space with its adjoint space admit the structure of a super operator system when equipped with the completely antiisometric involution given by antilinear extension of 
$\, v \otimes \overline w \mapsto w \otimes \overline v\, $. 
The projective tensor product $\, ( \mathcal V \widehat\otimes \overline{\mathcal V}\, ,\, 
1_{\mathcal V} \widehat\otimes 1_{\overline{\mathcal V}} )\, $ admits the structure of a pointed super operator space with the corresponding involution.
The unital projective tensor product (projective tensor product)
$\, \mathfrak X \widehat{\otimes }_u \mathfrak Y\, $ ($\, \mathfrak X \widehat{\otimes } \mathfrak Y\, $) and the injective tensor product 
$\, \mathfrak X \check\otimes\mathfrak Y\, $ of two super operator systems admit the structure  of a super operator system (pointed super operator space) equipped with the completely antiisometric involution given by antilinear extension of 
$\, x \otimes y \mapsto x^* \otimes y^*\, $.
\par\bigskip\noindent
{\it Proof.}\quad Considering the projective tensor product 
$\, \mathcal V \widehat\otimes \overline{\mathcal V}\, $, the assignment $\, v \otimes \overline w \mapsto w \otimes \overline v\, $ extends to a completely antilinear and antiisometric involution on the tensor product. Indeed, if 
$\, x = \alpha ( v \otimes \overline w ) \beta \in M_n ( \mathcal V \widehat\otimes \overline{\mathcal V} )\, $ with 
$\, v \in M_p ( \mathcal V )\, ,\, \overline w \in M_q ( \overline{\mathcal V} )\, ,\, \alpha \in M_{n , p \times q} ( \mathbb C )\, ,\, 
\beta \in M_{p \times q , n} ( \mathbb C )\, $, then 
$\, x^* = \overline\beta ( w \otimes \overline v ) \overline\alpha \, $ and 
$\, \Vert x^* \Vert \leq \Vert x \Vert\, $ from the definition of the projective norm. By symmetry one gets 
$\, \Vert x^* \Vert = \Vert x \Vert\, $. So the projective tensor product is naturally a pointed super operator space. Now the unification of a pointed super operator space 
$\, ( \mathfrak V\, ,\, *_{\mathfrak V} )\, $ is a super operator system.
Let 
$\,\rho : \mathfrak V  \rightarrow \mathcal B ( \mathcal H )\, $ be any unital representation, dropping to a representation of the unification. Then the composition of the involution on $\, \mathfrak V\, $ with the adjoint map of the algebra of bounded operators on Hilbert space defines another completely contractive unital representation which therefore drops to the unification. Since this is true for all such representations one finds that the given involution drops to a completely antiisometric unital involution on the unification of $\, \mathfrak V\, $ which is therefore a super operator system. 
In the same way one verifies that the second involution defined on the projective tensor product 
$\, \mathfrak X \widehat\otimes \mathfrak Y\, $ is completely antiisometric and drops to a completely antiisometric antilinear involution on the unital projective tensor product. The proof that the first of the outlined involutions is completely antiisometric for the Haagerup tensor product is very similar. It can also be derived using injectivity of the Haagerup tensor product on embedding the unital operator space $\, \mathcal V\, $ into an enveloping $C^*$-algebra $\, A\, $ and using the description of $\, A {\otimes }_h A\, $ as a subspace of the free product of two copies $\, A_i\, ,\, i = 1 , 2\, $ of $\, A\, $ amalgamated over the identity by linear extension of the assignment 
$\, a \otimes b \mapsto a_1 * b_2\, $ due to Christensen-Effros-Pisier-Sinclair (see \cite{Pa}, Theorem 17.6.). Then the involution of first type corresponds to composition of the $C^*$-involution on 
$\, A_1 * A_2\, $ with the order two automorphism exchanging the two copies of $\, A\, $, the first being completely antiisometric, and the second a completely isometric map. In case of the (spatial) injective tensor product the proof is obvious  \qed 
\par\bigskip\noindent
Besides the unital projective tensor product we wish to define another tensor product which is closely related with the Haagerup tensor norm. Unlike the former the Haagerup tensor product is not commutative so that the product type involution on $\,\mathfrak X \otimes \mathfrak Y\, $ is not completely antiisometric for the Haagerup tensor norm $\,\{ \Vert\cdot {\Vert }_{h , n} {\}}_{n\in\mathbb N}\, $. Let 
$\, x\mapsto x^*\, $ denote the product type involution on the algebraic tensor product (extended to matrices by composing it with the transpose mapping). Then one defines a symmetrized Haagerup tensor norm putting 
$$ \Vert x {\Vert }_{h^* , n}\> :=\> \max\, \{ \Vert x {\Vert }_{h , n}\, ,\, \Vert x^* {\Vert }_{h , n} \}\> . $$
It is easy to verify the axioms of a matrix norm. The corresponding completion of the algebraic tensor product is denoted $\, \mathfrak X {\otimes }_{h^*} \mathfrak Y\, $. One obtains a completely $*$-isometric unital representation of the resulting pointed super operator space by taking the direct sum of a completely isometric unital representation for the (unsymmetrized) Haagerup norm and the completely contractive representation obtained by composing the first representation with the involution on the left and the adjoint operation on the right with grading operator $\,\epsilon\, $ exchanging the two copies of the same underlying Hilbert space, showing that the symmetrized Haagerup tensor product is a super operator system. Although the symmetrized Haagerup tensor product is commutative, by the process of symmetrization associativity of the Haagerup tensor product is lost, i.e. 
$\, ( {\mathfrak V}\, {\otimes }_{h^*}\, {\mathfrak W} )\, {\otimes }_{h^*}\, {\mathfrak X}\, $ is in general different from $\, {\mathfrak V}\, {\otimes }_{h^*}\, ( {\mathfrak W}\, {\otimes }_{h^*}\, {\mathfrak X} )\, $. 
\par\smallskip\noindent
In the progress of our investigation of operator tensor products we shall need a version of the Christensen-Effros-Sinclair-Theorem suited for our purposes. Let $\,  A\, $ be a $C^*$-algebra and 
$\,\mathfrak V\, $ a super operator space which is also an 
(left) $A$-module, which is to say that there exists a completely bounded left action 
$\,\lambda : A \times \mathfrak V \rightarrow \mathfrak V\, $  meaning that the natural extensions 
$\, id_{M_n} \otimes\lambda  :  A \times M_n ( \mathfrak V ) \rightarrow M_n ( \mathfrak V )\, $ induced by 
$\,\lambda\, $ are uniformly bounded for all $\, n \in \mathbb N\, $. We only consider actions which are nondegenerate in the following sense:  given an approximate unit $\, \{ u_{\lambda } \}_{\lambda }\, $ for $\, A\, $, the net $\, \{ u_{\lambda } \cdot v \}_{\lambda }\, $ converges to 
$\, v\, $ for any $\, v\in \mathfrak V\, $. It follows from the nondegeneracy condition that the action extends naturally to a completely bounded action of the $C^*$-unitization $\,\widetilde A\, $ by setting 
$$ ( a , \lambda{\bf 1} )\cdot v = a\cdot v + \lambda v $$
using associativity of the action. Therefore in the following we may assume $\, A\, $ to be unital without loss of generality.
Putting 
$\, v \cdot a := ( a^*\cdot v^* )^*\, $ defines a right action of $\, A\, $ on 
$\, \mathfrak V\, $. We assume that the adjoint right action of $\, A\, $ is compatible with the left action in the sense that 
$\, a\cdot ( v \cdot a' ) = ( a\cdot x ) \cdot a'\, $ (if $\,\mathfrak V\, $ is a pointed super operator space resp. super operator system we also require that $\, a\cdot *_{\mathfrak V} = *_{\mathfrak V}\cdot a\, $).  
It then follows directly from the definition of the adjoint actions that they are again completely bounded. 
In our applications we shall usually need stronger requirements -- that the action is completely contractive in a certain sense. A condition which looks a bit farfetched at first but proves useful in a technical surrounding is the following: given a completely bounded action of $\, A\, $ on 
$\, \mathfrak V\, $ as above the action is said to be {\it u-weakly completely contractive} if for every unitary $\, u\in A\, $ the assignment 
$$ v \mapsto u\cdot v $$
is completely contractive (= completely isometric). It is said to be {\it weakly completely contractive} if for every 
$\, a\in A\, $ with $\, \Vert a \Vert \leq 1\, $ the assignment 
$$ v\mapsto a\cdot v $$
is completely contractive. The significance of the former notion is most easily recognized in its proper context. Namely let $\, \mathcal U\, $ denote a topological group (for example the unitary group of a 
$C^*$-algebra). A (pointed) super operator space $\, \mathfrak V\, $ is said to be a weak $\mathcal U$-bimodule, iff 
there is defined a left action 
$$ \mathcal U \times \mathfrak V\>\longrightarrow \mathfrak V\quad ,\quad 
\left( u\, ,\, v \right)\>\mapsto\> u\cdot v $$
which is continuous, linear and completely isometric in $\,\mathfrak V\, $, and associative, i.e. 
$\, u_1\cdot ( u_2\cdot v ) = (u_1u_2)\cdot v\, $, and such that putting 
$$ v\cdot u\> =\> ( u^{-1}\cdot v^* )^* $$
defines a right action of $\,\mathcal U\, $ on $\,\mathfrak V\, $ commuting with the left action 
(and if $\,\mathfrak V\, $ is pointed $\, u\cdot *_{\mathfrak V} = *_{\mathfrak V}\cdot u\, $ for any 
$\, u\in\mathcal U\, $).
The actions extend via matrix multiplication to completely bounded left and right actions of matrix algebras $\, M_n ( A )\, $ (resp. $\, {\mathcal U}^n\, $) on $\, M_n ( \mathfrak V )\, $ for each $\, n\, $. If the so defined extended actions are contractive (i.e. if the natural induced maps into 
$\, \mathcal C\mathcal B ( M_n ( \mathfrak V )\, ,\, M_n ( \mathfrak V ) )\, $ is contractive for each $\, n\, $)  the actions will be called {\it strongly completely contractive} and $\, \mathfrak V\, $ will be called a strong 
$A$-bimodule (resp. strong $\mathcal U$-bimodule), otherwise we say that $\, \mathfrak V\, $ is a (weak resp. u-weak) $A$-bimodule. For any super operator system $\, {\mathfrak X}_2\, $ 
(or super operator space) a {\it matrix decomposition} is synonymous for  a strongly completely contractive action of $\, \mathbb C\oplus \mathbb C\, $. It will be called an {\it $F$-matrix decomposition} if the strong action of 
$\, \mathbb C\oplus \mathbb C\, $ viewed as diagonal of $\, M_2 ( \mathbb C )\, $ with orthogonal projections $\, e\, $ and $\, 1-e\, $ corresponding to the left upper resp. lower right corner extends to an action of $\, M_2 ( \mathbb C )\, $ such that 
$$ F\> =\> \begin{pmatrix} 0 & 1 \\ 1 & 0 \end{pmatrix} $$
acts completely isometric, {\it and} such that the element 
$\, F = F\cdot {\bf 1}\, $ of the injective envelope 
$\, I ( {\mathfrak X}_2 )\, $ is a selfadjoint (even) unitary, {\it and} such that putting 
$\,\mathfrak X = (1-e)\cdot {\mathfrak X}_2  \cdot (1-e)\, $ the natural map 
$$ {\mathfrak X}_2 \rightarrow M_2 ( \mathfrak X ) $$
is completely contractive. 
To make things just a little more complicated we will also have to consider a notion which is slightly weaker than an $F$-matrix decomposition. A {\it weak $F$-matrix decomposition} of $\, {\mathfrak X}_2\, $ means a matrix decomposition  which extends to a completely bounded action of 
$\, M_2 ( \mathbb C )\, $ such that $\, F\, $ as above acts completely isometric restricted to any corner of the decomposition, {\it and} the element $\, F = F\cdot {\bf 1}\, $ is an even selfadjoint unitary of 
$\, I ( {\mathfrak X}_2 )\, $. 
\par\smallskip\noindent
The following definition of an {\it admissible action} is designed to accommodate the main applications that we have in mind, in order to get the strongest possible results. Some of the conditions are of a seemingly very special nature, and it may be that in a different situation another perhaps more simple-minded set of assumptions is sufficient. For example, if instead of a (weak) $F$-matrix decomposition one assumes the structure of a strong $ M_2 ( \mathbb C )$-bimodule the proof of the following 
SCES-Theorem will be much shorter and more direct but this does not fit with the situation encountered in our applications.  If $\,\mathfrak X\, $ is a super operator system a bimodule action of 
$\, A\, $ (resp. $\,\mathcal U\, $) on $\,\mathfrak X\, $ will be called {\it admissible} if 
it is u-weakly completely contractive and if the following conditions hold: for every unitary 
$\, u\in A\, $ (resp. $\, u\in\mathcal U\, $), putting $\, U = u \oplus 1\in A \oplus A\, $, there exists a super operator system $\, \overline{\mathfrak X}_{2 , U}\, $  and a unital complete $*$-contraction 
$\, p_U : \overline{\mathfrak X}_{2 , U} \twoheadrightarrow M_2 ( \mathfrak X )\, $ with dense image 
(it is in fact sufficient that $\, p_U\, $ is completely contractive when restricted to any enlarged corner 
$\, {\Delta }^{i j}_{2 , U}\, $ as defined below) such that the natural $F$-matrix decomposition of 
$\, M_2 ( \mathfrak X )\, $ lifts to a weak $F$-matrix decompositon on 
$\, \overline{\mathfrak X}_{2 , U}\, $, and there exists a subsystem 
$$ \overline{\mathfrak Y}_{2 , U}\> =\> \begin{pmatrix} \overline{\mathfrak Y} & F\cdot \overline{\mathfrak X} \\ \overline{\mathfrak X}\cdot F & \overline{\mathfrak X} \end{pmatrix}\> ,\quad \overline{\mathfrak Y}\subseteq F\cdot \overline{\mathfrak X}\cdot F $$
with $\,\overline{\mathfrak X} = ( 1 - e ) \overline{\mathfrak X}_{2 , U} ( 1 - e )\, $, 
such that the adjoint action of the unitary $\, U\, $ (with $\, A\oplus A\, $ viewed as subalgebra of diagonal elements in $\, M_2 ( A )\, $) on $\, p_U ( \overline{\mathfrak Y}_{2 , U} )\, $ given by 
$\, x\mapsto ad_U ( x ) =  U\cdot x\cdot U^*\, $ lifts to an action on $\, \overline{\mathfrak Y}_{2 , U}\, $
which is compatible with the matrix decomposition and trivial on the lower right corner 
$\, \overline{\mathfrak X}\, $. Putting 
$\, F_{U^n} = ad_{U^n} ( F )\, $ it is also assumed that the elements in the set 
$\, \{\, F\cdot F_{U^n}\, \vert\, n\in \mathbb Z\, \}\, $ are evenly graded and $ad_U$-invariant and that 
$\, F\cdot F_{U^n} = F\circ F_{U^n}\, $ (where $\circ\, $ denotes the $C^*$-product operation in the injective envelope). Let $\, {\Delta }^{ij}_{2 , U}\, $ denote the corner determined by the $i-th$ row and 
$j-th$ column of the matrix decomposition but with elements $\,\{ F_{U^n}\, ,\, F\cdot F_{U^n} \}\, $ adjoined, so that 
$$ F\cdot {\Delta }^{i j}_{2 , U} = {\Delta }^{(i+1) j}_{2 , U}\> ,\quad {\Delta }^{i j}_{2 , U}\cdot F = 
{\Delta }^{i (j+1)}_{2 , U} $$
with superscripts mod $2$. Also put $\, {\Delta }_{2 , U} = {\Delta }^{1 1}_{2 , U} + 
{\Delta }^{2 2}_{2 , U}\, $ and $\, {\nabla }_{2 , U} = F\cdot {\Delta }_{2 , U} = {\Delta }_{2 , U}\cdot F\, $.
The action of $\, F\, $ is still supposed to be completely isometric intertwining the larger subsystems 
$\, {\Delta }_{2 , U}\, $ and $\, {\nabla }_{2 , U}\, $ (which is consistent with our examples although probably some weaker hypotheses will also do). The map $\, ad_U\, $ is supposed to be completely isometric when restricted to the intersection of 
$\, \overline{\mathfrak Y}_{2 , U}\, $ with either of these subsystems (but not necessarily on the whole space). It follows from the condition 
$\, F\cdot F_{U^n} = F\circ F_{U^n}\, $ etc. that the natural map 
$$ \mu  : {\overline{\mathfrak X}}_{2 , U} \twoheadrightarrow M_2 ( \overline{\mathfrak X} ) $$
by identification of the corners via the $F$-action is completely isometric when restricted to any of the subspaces $\, {\Delta }^{ij}_{2 , U}\, $. The injective envelope of 
$\, \overline{\mathfrak Y}_{2 , U} \cap {\Delta }_{2 , U}\, $ (or more precisely the subspace of the injective envelope generated linearly by $C^*$-products from elements of 
$\,\overline{\mathfrak Y}_{2 , U} \cap {\Delta }_{2 , U}\, $ with the selfadjoint unitaries 
$\, \{ F_{U^n} \}\, $ from the left and right side) should contain 
$\, {\Delta }_{2 , U}\, $ by some completely isometric $*$-linear extension of the inclusion 
$\, \overline{\mathfrak Y}_{2 , U}\cap {\Delta }_{2 , U} \subseteq I ( \overline{\mathfrak Y}_{2 , U}\cap {\Delta }_{2 , U} )\, $. 
In addition there should be two unital complete contractions 
$$ {\nu }_l : {\overline{\mathfrak XÊ}}_{2 , U}\, \twoheadrightarrow\, {\overline{\mathcal V}}_{2 , U , l}\> ,\quad {\nu }_r : {\overline{\mathfrak X}}_{2 , U}\, \twoheadrightarrow\, {\overline{\mathcal V}}_{2 , U , r} $$
onto unital operator spaces with compatible (left and right) weak $F$-matrix decomposition such that in the first case the {\it left } action of $\, F\, $ is given by $C^*$-multiplication with $\, F\, $ from the left in the injective envelope (hence is a strong action), and in the second case the {\it right} action of $\, F\, $ is given by $C^*$-multiplication from the right with $\, F\, $
(in the remaining cases one again assumes that the action of $\, F\, $ is completely isometric restricted to  the images of $\, {\Delta }_{2 , U}\, $ and $\, {\nabla }_{2 , U}\, $, such that if 
$\, x\in M_n ( {\overline{\mathfrak X}}_{2 , U} )\, $ is contained in a subspace 
$\, {\Delta }^{i j}_{2 , U}\, $ one has 
$$ \Vert x {\Vert }_n\> =\> \max\,\{ \Vert {\nu }_l ( x ) {\Vert }_n\, ,\> \Vert {\nu }_r ( x ) {\Vert }_n \}\> . $$
The map $\, ad_U\, $ is supposed to drop to a completely isometric map on the images of the intersections of 
$\,\overline{\mathfrak Y}_{2 , U}\, $ with the (enlarged) diagonal/off-diagonal subspaces in both 
$\,  {\overline{\mathcal V}}_{2 , U , l}\, $ and $\,  {\overline{\mathcal V}}_{2 , U , r}\, $
(in particular the images of the elements $\,\{ F_{U^n} \}\, $ still define selfadjoint unitaries in the injective envelopes of both spaces which carry a unique structure of a unital $C^*$-algebra). 
\par\smallskip\noindent
The action of $\, A\, $ (or $\,\mathcal U\, $) on $\,\mathfrak X\, $ will be called {\it strongly admissible} if it is u-weakly completely contractive and if there exists a super operator system 
$\, \overline{\mathfrak X}\, $ which is a strong 
$A$- (or $\mathcal U$-)bimodule and a unital  $*$-linear complete contraction $\, p : \overline{\mathfrak X}Ê\twoheadrightarrow \mathfrak X\, $ with dense image commuting with the corresponding bimodule structures.
In particular a strong $A$-bimodule 
$\,\mathfrak X\, $ is strongly admissible taking $\,\overline{\mathfrak X} = \mathfrak X\, $. 
\par\bigskip\noindent
{\it Example.}\quad Let $\, \mathfrak V\, $ and $\, \mathfrak W\, $ be super operator spaces with 
$\, \mathfrak V\, $ a weak $A$-module and consider the projective tensor product 
$\, \mathfrak P = \mathfrak V \widehat{\otimes } \mathfrak W\, $. From the proof of Lemma 4 
$\,\mathfrak P\, $ can be given the structure of a super operator space with product involution of the 
involutions of $\, \mathfrak V\, $ and $\, \mathfrak W\, $ respectively. Consider the natural left action of 
$\, A\, $ on the algebraic tensor product
$\, \mathfrak V \otimes \mathfrak W\, $ equipped with the projective norm. We wish to show that this action is weakly completely contractive so that it extends to a weakly completely contractive action on the projective tensor product. Let $\, a \in A\, $ and $\, x \in M_n ( \mathfrak V \otimes \mathfrak W )\, $ be given and suppose that 
$\, x\, $ has a product decomposition of the form $\, \alpha\cdot ( v \otimes w )\cdot\beta\, $ with 
$\, \alpha\in M_{n , p\times q} ( \mathbb C )\, ,\, \beta\in M_{p\times q , n} ( \mathbb C )\, $ and 
$\, v\in M_p ( \mathfrak V )\, ,\, w\in M_q ( \mathfrak W )\, $. Then the element 
$\, a\cdot x\, $ has a decomposition 
$\, \alpha\cdot ( ( ( 1_p \otimes a )\cdot v ) \otimes w )\cdot\beta\, $ so that the projective norm of 
$\, a\cdot x\, $ is smaller than $\,  \Vert \alpha \Vert\cdot \Vert ( 1_p \otimes a )\cdot v\Vert\cdot 
\Vert w\Vert\cdot \Vert \beta \Vert \leq 
\Vert a \Vert\cdot \Vert\alpha\Vert\cdot\Vert v\Vert\cdot\Vert w\Vert\cdot\Vert\beta\Vert\, $ and taking the infimum over all such decompositions one gets 
$\, \Vert a\cdot x\Vert \leq \Vert a\Vert\cdot \Vert x\Vert\, $ as desired. By taking the adjoint right action one checks the compatibility conditions so that $\, \mathfrak P\, $ is a weak $A$-bimodule. If 
$\, \mathfrak W\, $ also happens to be a weak $B$-bimodule, then $\,\mathfrak P\, $ is a weak $A$- and 
weak $B$-bimodule, and each of the four (left and right) actions is compatible with each of the others.
\par\bigskip\noindent
{\bf Lemma 5.}\quad Suppose that the super operator space $\,\mathfrak V\, $ is a  u-weak $A$-bimodule. Then the adjoint action of the unitary group of $\, A\, $ extends to a  completely isometric and grading preserving action on its universal enveloping graded pre-operator system 
$\,\widehat{\mathfrak V}\, $. 
\par\bigskip\noindent
{\it Proof.}\quad The statement of the Lemma should be clear from the fact that given any completely contractive $*$-linear representation $\, \rho\, $ of $\,\mathfrak V\, $ the composition 
$\, {\rho }_u = \rho\circ ad_u\, $ for $\, u\in A\, $ unitary defines another $*$-linear representation\qed
 \par\bigskip\noindent
{\bf Theorem.}\, (Super-Christensen-Effros-Sinclair-Theorem)\quad
Let $\, A\, $ be a $C^*$-algebra (resp. $\,\mathcal U\, $ a topological group) and $\, \mathfrak X\, $ a super operator system which is a strongly 
admissible $A$-(resp. $\mathcal U$-)bimodule. Then the bimodule action extends to a graded action on its enveloping graded operator system $\,\widehat{\mathfrak X}\, $ which is strongly completely contractive. 
In particular $\, \mathfrak X\, $ and $\, {\mathfrak X}_0\, $ are strong $A$-(resp. $\mathcal U$-)bimodules and there exists a graded Hilbert space $\,\widehat{\mathcal H}\, $, a unital graded completely isometric $*$-representation 
$\, \pi : \widehat{\mathfrak X} \rightarrow\mathcal B ( \widehat{\mathcal H} )\, $ and a 
$C^*$-representation 
$\, \lambda : A \rightarrow \mathcal B ( \mathcal H )^{ev}\, $ 
(resp. unitary representation 
$\, \mu : \mathcal U\rightarrow \mathcal B ( \mathcal H )^{ev}\, $) 
into the subalgebra of elements of even degree such that 
$$ \pi ( a\cdot x\cdot a' )\> =\> \lambda ( a )  \pi ( x )  \lambda ( a' ) \> $$
or
$$ \pi ( u\cdot x\cdot u' )\> =\> \mu ( u ) \pi ( x ) \mu ( u')\>  $$
respectively.
\par\smallskip\noindent
If $\, \mathfrak X\, $ is an operator system which is an admissible $A$-bimodule 
($\mathcal U$-bimodule), then $\,\mathfrak X\, $ is a strong $A$-(resp. $\mathcal U$-)bimodule (operator $A$-system) and there exist compatible representations $\,\pi\, $ of $\,\mathfrak X\, $ and 
$\,\lambda\, $ of $\, A\, $ (resp. $\,\mu\, $ of $\,\mathcal U\, $)  on an (ungraded) Hilbert space as above.
\par\bigskip\noindent 
{\it Proof.}\quad We only write out the proof in case of the action by a $C^*$-algebra $\, A\, $ since the argument in case of a toplological group $\,\mathcal U\, $ is completely analogous. We begin with the second part so let $\, \mathfrak X\, $ be an operator system which is an admissible $A$-bimodule. Let 
$\, u\in A\, $ be a unitary and put 
$$ U\> =\> \begin{pmatrix} u & 0 \\ 0 & 1 \end{pmatrix}\>\in M_2 ( A )\> . $$
Then there exists a super operator system 
$\, \overline{\mathfrak X}_{2 , U}\, $ and a unital complete $*$-contraction with dense image 
$\, p_U : \overline{\mathfrak X}_{2 , U} \longrightarrow M_2 ( \mathfrak X )\, $
such that the adjoint action of $\, U\, $ on $\, M_2 ( \mathfrak X )\, $ lifts to a unital complete $*$-isometry on $\, \overline{\mathfrak Y}_{2 , U}\subseteq \overline{\mathfrak X}_{2 , U}\, $ (or its intersection with the subsystems $\, {\Delta }_{2 , U}\, ,\, {\nabla }_{2 , U}\, $) and admitting a compatible weak $F$-matrix decomposition. 
Any completely positive extension $\, \phi : I ( {\mathfrak X}_{2 , U} ) \rightarrow M_2 ( I ( \mathfrak X ) )\, $ 
of $\, p_U\, $ to the injective envelopes will be a $\mathbb C \oplus \mathbb C$-module map, since the strong action of $\,\mathbb C\oplus\mathbb C\, $ extends uniquely by rigidity to a strong action on the injective envelopes where it is implemented by $C^*$-multiplication with elements from the corresponding subalgebras which are identified by the completely positive map $\,\phi\, $, so that in any case $\,\phi\, $ is compatible with the matrix decomposition. By the same argument any extension 
$\,\phi\, $ is necessarily a $C^* ( F_{U^n} )$-module map for each $\, n\in\mathbb Z\, $.
We wish to show that $\, ad_U\, $ acts completely positive on $\, M_2 ( \mathfrak X )\, $. Since the adjoint action of $\, U\, $ on $\, \overline{\mathfrak Y}_{2 , U}\, $ is unital and completely isometric it follows that 
$\, ad_U\, $ extends uniquely to a graded $C^*$-automorphism on 
$\, I ( \overline{\mathfrak X}_{2 , U} ) = I ( \overline{\mathfrak Y}_{2 , U} )\, $ as well as 
$\, I ( {\Delta }_{2 , U} )\, $ and  $\, I ( {\nabla }_{2 , U} )\, $ by rigidity and Proposition 1. 
On the other hand the left action of $\, u\, $ is completely isometric on $\, \mathfrak X\, $ so the map 
$\, ad_U\, $ is unital and completely isometric restricted to the subsystem 
$$ {\mathcal S}_{\mathfrak X}\, =\, \begin{pmatrix} 
\mathbb C {\bf 1} & \mathfrak X \\ \mathfrak X & \mathbb C {\bf 1} \end{pmatrix} \subseteq 
M_2 ( \mathfrak X ) \> . $$
Then it extends to a unique $C^*$-automorphism on the injective envelope 
$\, I ( {\mathcal S}_{\mathfrak X} )\, $ which admits a matrix decomposition of the form 
$$ \begin{pmatrix} I_{11} & I ( {\mathfrak X} ) \\ I ( {\mathfrak X} ) & I_{22} \end{pmatrix}\> . $$
We claim that it is isomorphic to $\, M_2 ( I ( {\mathfrak X} ) )\, $.
Indeed it contains the element 
$$ \begin{pmatrix} 0 & 1 \\ 1 & 0 \end{pmatrix} $$
whose square is the identity element. Then define a product on the subspace $\, I ( {\mathfrak X} )\, $ sitting in the right upper corner of the injective envelope of $\, {\mathcal S}_{\mathfrak X}\, $ by the rule 
$$ \begin{pmatrix} 0 & x \\ 0 & 0 \end{pmatrix} * \begin{pmatrix} 0 & y \\ 0 & 0 \end{pmatrix} \> :=\> 
\begin{pmatrix} 0 & x \\  0 & 0 \end{pmatrix}\circ \begin{pmatrix} 0 & 1 \\ 1 & 0 \end{pmatrix}\circ 
\begin{pmatrix} 0 & y \\ 0 & 0 \end{pmatrix} $$
where $\,\circ\, $ denotes the $C^*$-multiplication in $\, I ( {\mathcal S}_{\mathfrak X} )\, $.
It is clear from the definition that the new product is associative making $\, I ( \mathfrak X )\, $ into a unital algebra, in fact a unital operator algebra.
Being an injective operator space it is equal to its own injective envelope. Since it is also unital it has a unique structure of $C^*$-algebra which is compatible with, hence equal to, its operator algebra structure. Then one sees that 
$\, I ( {\mathcal S}_{\mathfrak X} )\, $ contains (and is equal to) the subspace 
$$ \begin{pmatrix} 0 & I ( {\mathfrak X} ) \\ I ( \mathfrak X ) & 0 \end{pmatrix} 
\> +\> \begin{pmatrix} 0 & 1 \\ 1 & 0 \end{pmatrix}\circ
\begin{pmatrix} 0 & I ( \mathfrak X ) \\ I ( \mathfrak X ) & 0 \end{pmatrix} $$
which being a $C^*$-algebra must be isomorphic with $\, M_2 ( I ( \mathfrak X ) )\, $.
Therefore the map $\, ad_U\, $ on $\, {\mathcal S}_{\mathfrak X}\, $ has a unique ($C^*$-automorphism) extension to 
$\, I ( {\mathcal S}_{\mathfrak X} ) \simeq M_2 ( I ( \mathfrak X ) )\, $ which commutes with the matrix decomposition. It remains to show that this extension restricted to $\, M_2 ( \mathfrak X ) \subseteq 
M_2 ( I ( \mathfrak X ) )\, $ is equal to $\, ad_U\, $. 
The completely isometric left and right action of $\, F\, $ mapping $\, {\nabla }_{2 , U}\, $ onto 
$\, {\Delta }_{2 , U}\, $ extends by rigidity to completely isometric  bijections $\, {\overline l}_F\, $ and 
$\, {\overline r}_F\, $ of the injective envelopes, so that putting 
$$ {\overline\lambda }_F ( x )\> :=\> F\, {\overline l}_F ( x )\quad ,\quad {\overline\rho }_F ( x )\> :=\> 
{\overline r}_F ( x )\, F   $$
defines, being unital bijective and completely isometric $C^* ( F )$-module maps, (ungraded) $C^*$-isomorphisms of 
$\, I ( {\nabla }_{2 , U} )\, $ and $\, I ( {\Delta }_{2 , U} )\, $. 
One also has the graded $C^*$-automorphisms
$$ Ad\, F ( x )\> =\> F\, x\, F\quad ,\quad ad_F ( x )\> =\> F\cdot x\cdot F $$
on the injective envelopes of both of the two subsystems as above.
Each of these is a $C^* ( F )$-module map and respects the $F$-matrix decomposition.
Let $\,\alpha\, $ denote the grading automorphism of $\, I ( \overline{\mathfrak X}_{2 , U} )\, $.
Put $\, {\underline{\mathfrak X}}_{2 , U} = \alpha ( {\overline{\mathfrak X}}_{2 , U} )\, $ and 
$\, \underline{\mathfrak X} = ( 1 - e )\, {\underline{\mathfrak X}}_{2 , U}\, ( 1 - e )\, $ with enveloping graded operator system $\, \widehat{\mathfrak X} = \overline{\mathfrak X} + \underline{\mathfrak X}\, $. Also put 
$\, {\underline l}_F = \alpha\circ {\overline l}_F\circ \alpha\, ,\, {\underline r}_F = \alpha\circÊ{\overline r}_F\circ\alpha\, $ etc. .
The map
$$  \overline\mu  :\>  \begin{pmatrix} ( {\overline l}_F\circ {\overline r}_F ) ( \widehat{\mathfrak X} ) & 
{\overline l}_F ( \widehat{\mathfrak X} ) \\ {\overline r}_F ( \widehat{\mathfrak X} ) & 
\widehat{\mathfrak X} \end{pmatrix}\> \twoheadrightarrow\> 
M_2 ( \widehat{\mathfrak X}  ) $$
by identification of the corners via the $F$-matrix decomposition is again completely isometric restricted to some enlarged corner.
The domain of $\, \overline\mu \, $ needn't be invariant under $\, ad_U\, $ but it is restricted to 
$\, {\overline{\mathfrak Y}}_{2 , U}\, $.
Of course the same arguments apply with respect to the $F$-matrix decomposition obtained from the maps 
$\, {\underline l}_F\, $ and $\, {\underline r}_F\, $ and replacing $\, {\overline{\mathfrak Y}}_{2 , U}\, $ by 
$\, {\underline{\mathfrak Y}}_{2 , U} = \alpha ( {\overline{\mathfrak Y}}_{2 , U} )\, $. 
The corresponding map to $\, M_2 ( \widehat{\mathfrak X} )\, $ is denoted $\, \underline\mu \, $.
The map 
$\, {\nu }_l \oplus {\nu }_r\, $ is completely isometric restricted to the enveloping operator system of 
any subspace $\, {\Delta }^{i j}_{2 , U}\, $ (since the unital complete isometry extends to a selfadjoint complete isometry of injective envelopes). Put 
$\, Ad_{F_{U^n}} ( x ) = F_{U^n} ad_{F_{U^n}} ( x ) F_{U^n}\, $.
Also put
$$ {\mathcal S}^{l , n}_{\widehat{\mathfrak X}}\> =\>  \begin{pmatrix} \mathbb C\, e & 
{\overline l}_{F_{U^n}} ( \widehat{\mathfrak X} ) \\  F_{U^n}\, {\overline l}_{F_{U^n}} 
( \widehat{\mathfrak X} )\, F_{U^n} & \mathbb C\, ( 1 - e ) \end{pmatrix}\> ,
\quad {\mathcal S}^{r , n}_{\widehat{\mathfrak X}}\> =\> ( \alpha\circ Ad_{F_{U^n}} ) ( {\mathcal S}^{l , n}_{\widehat{\mathfrak X}} ) \> . $$
By the general identity $\, {\overline\lambda }_F ( \alpha ( x ) ) = \alpha ( {\overline\rho }_F ( x ) )\, $ one gets 
$$  {\overline l}_{F_{U^n}} ( \alpha ( x ) )\> =\> ( \alpha\circ Ad\, F_{U^n}\circ {\overline r}_{F_{U^n}} ) ( x )\> =\> (Ê\alpha\circ Ad_{F_{U^n}}\circ {\overline l}_{F_{U^n}} ) ( x ) $$
so that $\, {\mathcal S}^{l / r , n}_{\widehat{\mathfrak X}}\, $ are operator systems with each of the (off diagonal) corners completely isometric to the operator system $\, \widehat{\mathfrak X}\, $. Their injective envelopes are each naturally isomorphic with $\, M_2 ( I ( \widehat{\mathfrak X} ) )\, $ by an argument as above and the same as the injective envelopes of the subspaces 
$$  {\mathcal U}_{\mathcal T}^{l , n}\> =\>  \begin{pmatrix} \mathbb C\, e 
& {\overline l}_{F_{U^n}} ( \widehat{\mathfrak X} ) \\  0 & \mathbb C\, ( 1 - e ) \end{pmatrix}\> , $$
and 
$$  {\mathcal L}_{\mathcal T}^{r , n}\> =\>  \begin{pmatrix} \mathbb C\, e & 0 
\\  {\underline r}_{F_{U^n}} ( \widehat{\mathfrak X} ) & \mathbb C\, ( 1 - e ) \end{pmatrix}\>  $$
of  upper resp. lower triangular matrices. We also consider the corresponding subspaces 
$$  {\mathcal U}_{\mathcal T , 0}^l \> =\>  \begin{pmatrix} \mathbb C\, e 
& {\overline l}_F ( \overline{\mathfrak X} ) \\  0 & \mathbb C\, ( 1 - e ) \end{pmatrix}\> , $$
and 
$$  {\mathcal L}_{\mathcal T , 0}^r\> =\>  \begin{pmatrix} \mathbb C\, e & 0 
\\  {\underline r}_F ( \underline{\mathfrak X} ) & \mathbb C\, ( 1 - e ) \end{pmatrix} $$
which are mutually adjoint to each other.
Then consider the unital operator spaces 
$$ \widetilde{\mathcal S}^n\> =\> {\overline\lambda }_{F_{U^n}} 
( {\mathcal S}^{l , n}_{\widehat{\mathfrak X}} )\>  =\> 
{\underline\rho }_{F_{U^n}} ( {\mathcal S}^{r , n}_{\widehat{\mathfrak X}} )\> ,\> 
\widetilde{\mathcal U}^{l , n}_{\mathcal T} = {\overline\lambda }_{F_{U^n}} ( {\mathcal U}^{l , n}_{\mathcal T} )\> ,\> \widetilde{\mathcal L}^{r , n}_{\mathcal T} = {\underline\rho }_{F_{U^n}} ( {\mathcal L}^{r , n}_{\mathcal T} ) $$
and their corresponding subspaces 
$$ \widetilde{\mathcal U}^{l , n}_{\mathcal T , 0}\> ,\quad 
\widetilde{\mathcal L}^{r , n}_{\mathcal T , 0} \> . $$
Choose a completely isometric realization of the injective envelope $\, I ( \widehat{\mathfrak X} )\, $ in 
$\, ( 1 - e )\, I ( {\overline{\mathfrak X}}_{2 , U} )\, ( 1 - e )\, $ which is elementwise fixed under $\, ad_U\, $. 
This is achieved starting from some arbitrary completely positive extension $\, j\, $ to injective envelopes of the $ad_U$-covariant embedding 
$$ \widehat{\mathfrak X}\>\hookrightarrow {\widehat{\overline{\mathfrak X}}}_{2 , U} $$ 
by "averaging" over the $ad_U$-translates of $\, j\, $ (compare with the discussion below).
Since 
$\, \widehat{\mathfrak X} \subseteq 
I ( \widetilde{\mathcal S}^0 )\, $ (by any unital completely positive embedding of the latter into 
$\, I ( {\overline{\mathfrak X}'}_{2 , U} )\, $ which must be a $C^* ( F )$-module map) the lower right corner of the latter is an injective envelope of the former and may be assumed to agree with the specified copy of $\, I ( \widehat{\mathfrak X} )\, $ as above. In the same way 
$\, \widehat{\mathfrak X} \subseteq I ( \widetilde{\mathcal S}^n )\, $
and one may choose both injective envelopes to contain the same copy of 
$\, I ( \widehat{\mathfrak X} )\, $ in the lower right corner. 
One checks that each element $\, F_{U^n}\, $ defines a selfadjoint unitary in the injective envelope of 
$\, {\mathcal U}^l_{\mathcal T , 0}\, $ and that any completely positive extension of the inclusion into 
$\, {\mathcal S}^{l , n}_{\widetilde{\mathfrak X}'}\, $ to injective envelopes is a $C^* ( F_{U^n} )$-module map. Then one gets the equality
$$  F_{U^n}\, F\> =\> F\, F_{U^{-n}}\> \in\>  
I ( \widehat{\mathfrak X} )Ê\oplus F\, I ( \widehat{\mathfrak X} )\, F  $$
by $ad_U$-invariance of any element in $\, I ( \widehat{\mathfrak X} )\, $ and also 
$$ \overline{\lambda }_F ( F_{U^n} )\> =\> F_{U^n}\> ,\> \overline{\rho }_F ( F_{U^n} )\> =\>
 F_{U^n}\quad \Rightarrow\quad Ad_F ( F_{U^n} )\> =\> F_{U^n} $$ 
for every $\, n\, $. It follows that (one may choose corresponding injective envelopes such that) 
$\, ( 1 - e )\, I ( \widetilde{\mathcal U}^{l , n}_{\mathcal T , 0} )\, ( 1 - e ) \subseteq I ( \widehat{\mathfrak X} )\, ,\, ( 1 - e )\, I ( \widetilde{\mathcal L}^{r , n}_{\mathcal T , 0} )\, ( 1 - e )\subseteq 
I ( \widehat{\mathfrak X} )\, $ which means that $\, ad_U\, $ leaves these subspaces elementwise fixed. 
This implies that 
$$ ( \overline{\lambda }_{F_{U^n}}\circ \overline{\lambda }_F ) ( ( 1 - e )\, I ( \widetilde{\mathcal U}^{l , 0}_{\mathcal T , 0} )\, ( 1 - e ) ) \subseteq 
I ( \widehat{\mathfrak X} ) $$
and since $\, \overline{\lambda }_{F_{U^n}} ( F_{U^k} ) = F_{U^k}\, $ for every $\, n\, $ and $\, k\, $ the same holds for all other corners, for example
$$ ( \overline{\lambda }_{F_{U^n}}\circ \overline{\lambda }_F ) ( F_{U^k}\, ( 1 - e )\, I ( \widetilde{\mathcal U}^{l , 0}_{\mathcal T , 0} )\, ( 1 - e ) ) \subseteq F_{U^k}\, I ( \widehat{\mathfrak X} ) \>  $$
etc. so every $\, I ( \widetilde{\mathcal U}^{l , n}_{\mathcal T , 0} )\, $ is contained in the injective envelope of each $\,  \widetilde{\mathcal S}^k\, $, and all of these agree with the injective envelope of 
$\, \widetilde{\mathcal S}^0\, $ which in particular follows to be $ad_U$-invariant as a subspace of 
$\, I ( {\overline{\mathfrak X}}_{2 , U} )\, $. 
Then for $\, x\in {\mathcal U}^l_{\mathcal T , 0}\, $ put 
$\, {\widetilde x}_n = \overline{\lambda }_{F_{U^n}} ( x )\in 
\widetilde{\mathcal U}^{l , n}_{\mathcal T , 0}\, $. 
We shall see that there exists a completely isometric copy of 
$\, {\mathcal U}^l_{\mathcal T , 0}\, $ such that its injective envelope admits an $ad_U$-covariant realization as a subspace of the direct sum of two copies of $\, I ( {\overline{\mathfrak X}}_{2 , U} )\, $  (with $\, ad_U\, $ acting diagonally)  and such that its lower right corner with respect to the matrix decomposition is contained in 
$\, I ( \widehat{\mathfrak X} ) \oplus I ( \widehat{\mathfrak X} )\, $ (in particular it is elementwise fixed under the action of $\, ad_U\, $). If 
$\, x\in {\mathcal U}^l_{\mathcal T , 0}\, $ is an element consider the sequence of elements 
$$ \left\{ {\overline x}_m\,\Bigm\vert\, 
{\overline x}_m\> =\> {1\over 2m+1}\,\sum_{k = -m}^m\, {\widetilde x}_k\,\right\} \> . $$ 
Let $\, \overline{\mathcal U}^{l , m}_{\mathcal T , 0}\, $ be the subspace generated by the elements 
$\, \{ {\overline x}_m\,\vert\, x\in {\mathcal U}^l_{\mathcal T , 0} \}\, $ so one gets a unital complete contraction 
$$ {\overline\sigma }_m :\> {\mathcal U}^l_{\mathcal T , 0}\> \twoheadrightarrow\> 
\overline{\mathcal U}^{l , m}_{\mathcal T , 0}\> ,\quad {\overline\sigma }_m ( x )\> =\> {\overline x}_m  $$
for each $\, m\in\mathbb N\, $. It is in fact completely bounded from below by the image of the map 
$\, {\nu }_l\, $. The idea is to construct a 
"halfsided" inverse to the map $\, \overline{\sigma }_m\, $ which when composed with the map 
$\, {\nu }_l\, $ will be completely contractive. 
The complete contraction 
$\, {\nu }_l : {\overline{\mathfrak X}}_{2 , U} \twoheadrightarrow  {\overline{\mathcal V}}_{2 , U ,l}\, $ has a completely contractive $ad_U$-invariant extension
$$ {\overline\psi }_l  : I ( {\overline{\mathfrak X}}_{2 , U} ) \> \longrightarrow\> 
I ( {\overline{\mathcal V}}_{2 , U , l}  ) $$
so that $\, {\psi }_l ( {\lambda }_{F_{U^n}} ( x ) ) = {\psi }_l ( x )\, $ for all $\, x\in {\overline{\mathfrak X}}_{2 , U}\, $ and $\, n\in\mathbb Z\, $ by the assumption on 
$\, {\overline{\mathcal V}}_{2 , U , l}\, $. 
This means that each $\, {\widetilde x}_k\, $ is sent  to the image of $\, x\, $ under $\, {\overline\psi }_l\, $, and the composition 
$\, {\overline\psi }_l\circ {\overline\sigma }_m\, $ equals $\, {\nu }_l\, $ on 
$\, {\mathcal U}^l_{\mathcal T , 0}\, $. In particular 
$\, \Vert {\overline\sigma }_m ( x ) {\Vert }_n \geq \Vert {\nu }_l ( x ) {\Vert }_n\, $. 
For $\, y = x^* \in {\mathcal L}^r_{\mathcal T , 0}\, $ one analogously constructs the sequence of elements 
$$ \left\{\, {\underline y}_m \,\Bigm\vert\, {\underline y}_m\> =\> {1\over 2m + 1}\,\sum_{-m}^m\, {\widetilde y}_m\, \right\} \> ,\quad {\widetilde y}_m\> =\> {\underline\rho }_{F_{U^n}} ( y )\> . $$
This amounts to a complete contraction 
$$ {\underline\sigma }_m :\> {\mathcal U}^l_{\mathcal T , 0}\>\twoheadrightarrow \> 
( \underline{\mathcal L}^{r , m}_{\mathcal T , 0} )^*\> ,\quad {\underline\sigma }_m ( x )\> =\> 
{\underline y}_m^* \> . $$
Using a similar argument gives  
$\, \Vert {\underline\sigma }_m ( x ) {\Vert }_n \geq 
\Vert {\nu }_r ( x ) {\Vert }_n\, $.
One concludes that $\, {\sigma }_m = {\overline\sigma }_m \oplus {\underline\sigma }_m\, $ is completely isometric for each $\, m\, $.
The sequence of maps $\, \{ {\sigma }_m \}\, $ must have an adherence value in the {\it bounded weak topology} on the set of completely contractive maps of $\, {\mathcal U}^l_{\mathcal T , 0}\, $ into 
$\, I ( \widetilde{\mathcal S}^0\oplus \widetilde{\mathcal S}^0 ) \simeq M_2 ( I ( \widehat{\mathfrak X} ))\oplus M_2 ( I ( {\widehat{\mathfrak X}}) )\, $ (induced by the so called {\it bounded weak topology} of completely contractive linear maps into $\, \mathcal B ( \mathcal H )\, $ for a given representation space of $\, M_2 ( I ( \widehat{\mathfrak X} ) )\, $, compare with chap. 7 of \cite{Pa}, employing some $ad_U$-invariant completely positive projection of the former with range equal to the latter space; to construct it one first extends the map $\, ad_U\, $ to a completely contractive map on $\, \mathcal B ( \mathcal H )\, $ by injectivity and then uses a simple averaging process starting with an arbitrary completely positive projection $\, p\, $ onto $\, M_2 ( I ( \widehat{\mathfrak X} ) )\, $ putting 
$\, {\overline p}_m  = \sum_{k=1}^m\, ( ad_{U^{-k}}\circ p\circ ( ad_U )^k )\, $ to obtain a sequence of approximately $ad_U$-invariant projections as $\, m\to\infty\, $ which again has some $ad_U$-invariant  limit point in the bounded weak topology) , which subset is compact, and choosing a convergent subsequence $\,\{ {\sigma }_{m_k}\,\vert\, k\in\mathbb N \}\, $, the limit map 
$$ \sigma\> =\> \lim_{k\to\infty } {\sigma }_{n_k} $$
is again a unital completely isometric  and the image space
$\, \sigma ( {\mathcal U}^l_{\mathcal T , 0} )\, $
is $ad_U$-invariant, since for each element $\, x\, $
$$ \Vert {\sigma }_m ( ad_U ( x ) ) - ad_U ( {\sigma }_m ( x ) ) \Vert\> \rightarrow\> 0\> , \quad m\to\infty\> .  $$
The injective envelope of $\, \sigma ( {\mathcal U}^l_{\mathcal T , 0} )\, $ admits an $ad_U$-covariant  realization such that its lower right corner is contained in 
$\, I ( \widehat{\mathfrak X} ) \oplus I ( \widehat{\mathfrak X} )\, $ which means that each of its elements is fixed under $\, ad_U\, $, so the same must be true for the injective envelope of 
$\, {\mathcal U}^l_{\mathcal T , 0}\, $.
The composite map 
$$ {\mathcal U}^l_{\mathcal T , 0}\largerightarrow {\mathcal S}_{\overline{\mathfrak X}} \buildrel p_U \over\largerightarrow {\mathcal S}_{\mathfrak X} $$
extends by a simple averaging process as above to a completely positive map $\,\phi\, $ on injective envelopes which commutes with the maps $\, ad_U\, $. Then $\,\phi\, $ is necessarily a 
$C^* ( F_{U^n} )$-module map for every $\, n\, $, so that if $\, x\in\mathfrak X\, $ one finds (approximately) a lift for 
$$ \begin{pmatrix} 0 & 0 \\ 0 & x \end{pmatrix}  $$
of the form
$$ F\, \begin{pmatrix} 0 & \overline x \\ 0 & 0 \end{pmatrix}\>\in\> 
F\, {\mathcal U}^l_{\mathcal T , 0}   $$
($C^*$-multiplication in $\, I ( {\mathcal U}^l_{\mathcal T , 0} )\, $).
Since the lifted element is $ad_U$-invariant and 
$\,\phi\, $ commutes with $\, ad_U\, $, then also its image must be $ad_U$-invariant.
Therefore, by rigidity,  $\, ad_U\, $ is the identity restricted to the subalgebra 
$$ \begin{pmatrix} \mathbb C & 0 \\ 0 & I ( \mathfrak X )\end{pmatrix} $$
of $\, M_2 ( I ( \mathfrak X ) )\, $ and thus a 
$\mathbb C \oplus I ( \mathfrak X )$-module map.
Following the proof of Theorem 15.12. of \cite{Pa} one gets that the map 
$$ A \longrightarrow I ( \mathfrak X )\quad ,\quad a\mapsto a\cdot {\bf 1} $$
of $\, A\, $ into the injective envelope of $\, \mathfrak X\, $ is a $C^*$-homomorphism implementing the action of $\, A\, $. Then the action must be strongly completely contractive so that $\, \mathfrak X\, $ is an operator $A$-system which proves the second assertion of the theorem 
(choosing a unital $C^*$-representation of $\, I ( \mathfrak X )\, $). 
\par\smallskip\noindent
Next assume that $\,\mathfrak X\, $ is a super operator system with strongly admissible $A$-action and let $\,\overline{\mathfrak X}\, $ be the given super operator system with strong $A$-action and $*$-linear complete contraction 
$\, p : \overline{\mathfrak X} \twoheadrightarrow \mathfrak X\, $ with dense image commuting with the 
$A$-actions. 
Then $\, p\, $ extends naturally to a graded map of enveloping graded operator systems 
$$ \widehat p\> .\> \widehat{\overline{\mathfrak X}}\> \longrightarrow\> \widehat{\mathfrak X} $$
which also is a complete contraction with dense image. We begin by showing that 
the $A$-action extends to a graded strongly completely contractive action on
$\, \widehat{\overline{\mathfrak X}}\, $. 
Let $\, U\in M_n ( A )\, $ be a unitary. Then the map $\, Ad\, U\, $ is unital, $*$-linear and completely contractive on $\, M_n ( \overline{\mathfrak X} )\, $, hence it extends to a graded completely contractive map on $\, M_n ( \widehat{\overline{\mathfrak X}} )\, $. 
One then gets that the left (and hence also right) action of $\, A\,$ extends to a graded strongly completely contractive action on
$\, \widehat{\overline{\mathfrak X}}\, $. This follows since for $\, U\in M_n ( A )\, $ unitary putting 
$$ \overline U\> =\> \begin{pmatrix} U & 0 \\ 0 & 1 \end{pmatrix} $$
the map $\, Ad\, \overline U\, $ on $\, M_{2n} ( \overline{\mathfrak X} )\, $ extends (being unital, $*$-linear and completely contractive) to a $\mathbb Z_2$-graded completely contractive map on 
$\, M_{2n} ( \widehat{\overline{\mathfrak X}} )\, $. But then the left action of $\, U\, $ on 
$\, \widehat{\overline{\mathfrak X}}\, $ must also be graded and completely contractive. 
Then for any normal element $\, X\in M_n ( A )\, $ consider the unitary 
$$ U\quad =\quad \begin{pmatrix} \sqrt{ 1 - XX^*} & X \\ X^* & - \sqrt{ 1 - XX^*}Ê\end{pmatrix} \>\in\> 
M_{2n} ( A ) \> . $$
Since $\, U\, $ acts completely contractive on $\, M_{2n} ( \widehat{\overline{\mathfrak X}} )\, $ the element 
$\, X\, $ acts completely contractive on $\, M_n ( \widehat{\overline{\mathfrak X}} )\, $. This holds in particular for selfadjoint elements, and then one easily extends to arbitrary elements so  the action of 
$\, A\, $ on $\,\widehat{\overline{\mathfrak X}}\, $ is strongly completely contractive and 
$\mathbb Z_2$-graded.
One then deduces that the left action of $\, A\, $ on $\, \mathfrak X\, $ extends to a graded u-weakly completely contractive action on $\, \widehat{\mathfrak X}\, $. To see this note that the left action of a unitary $\, u\, $ extends (being completely isometric) uniquely to a completely isometric map $\, l_u\, $ on the injective envelope $\, I ( \mathfrak X )\, $ which is a unital $\mathbb Z_2$-graded $C^*$-algebra.
The corresponding map on $\, I ( \overline{\mathfrak X} )\, $ is certainly $\mathbb Z_2$-graded (in any case when restricted to $\, \widehat{\overline{\mathfrak X}}\, $). Let $\, \pi : I ( \overline{\mathfrak X} ) \rightarrow I ( \mathfrak X )\, $ be any unital completely contractive extension of $\, p\, $, and consider the composite map 
$$ I ( \overline{\mathfrak X} ) \buildrel l_{u^*}\over\longrightarrow I ( \overline{\mathfrak X} ) \buildrel \pi\over\longrightarrow I ( \mathfrak X ) \buildrel l_u\over\longrightarrow I ( \mathfrak X )\> . $$
It is another completely contractive extension of $\, p\, $, hence necessarily selfadjoint for the ordinary adjoint operation. Since it is also selfadjoint for the superinvolution (at least when restricted to 
$\,\overline{\mathfrak X}\, $) it must restrict to a $\mathbb Z_2$-graded map from 
$\,\widehat{\overline{\mathfrak X}}\, $ to $\,\widehat{\mathfrak X}\, $. But then, since the maps 
$\, l_{u^*}\, $ and $\,\pi\, $ are $\mathbb Z_2$-graded, so is the map $\, l_u : \widehat{\mathfrak X} \rightarrow \widehat{\mathfrak X}\, $. From this one infers that the left action of $\, A\, $ on 
$\, \mathfrak X\, $ extends to a graded and u-weakly completely contractive action on its enveloping graded operator system. Clearly this action is admissible so the result follows from the argument above for the ungraded case \qed 
\par\bigskip\noindent
The previous Theorem gives some information about the relation between the maximal $C^*$-tensor product $\, A\, {\otimes }_{max}\, B\, $, which is the closure of the algebraic tensor product in the direct sum over all joint commuting representations of $\, A\, $ and $\, B\, $, 
and the symmetrized Haagerup tensor product 
$\, A\, {\otimes }_{h^*} B\, $ of two given (unital) $C^*$-algebras $\, A\, $ and $\, B\, $. One can deduce that the former is completely $*$-isometric with the largest operator system $\, {\mathfrak X}_0\, $ which is smaller than the super operator system $\, A\, {\otimes }_{h^*} B \, $ and admits a u-weak $A$- and 
$B$-bimodule structure. This follows since such an operator system is automatically admissible for both actions, and hence is necessarily a strong 
$A\, {\otimes }_{max} B$-bimodule by the Theorem above. To see that it is admissible consider first the super operator system $\,\mathfrak Z = A\, {\otimes }_{h^*} B\, $. From the very definition of the symmetrized Haagerup tensor product it is easy to see that the natural left and right $A$-action is weakly completely contractive, this being the case for the (unsymmetrized) Haagerup tensor product (the case of the $B$-action follows by symmetry). 
A completely isometric $*$-representation of the symmetrized Haagerup tensor product is constructed from a unital completely isometric representation of the unsymmetrized Haagerup tensor product as indicated in the example above. On the other hand a representation of the ordinary Haagerup tensor product is given by embedding into the amalgamated free product 
$\, A\, {*}_{\mathbb C} B\, $ where an elementary tensor $\, a \otimes b\, $ is represented by the ordered product  $\, a * b\, $. Then for the symmetrized Haagerup tensor product the elementary tensor is represented by 
the element $\, ( a * b ) \widehat{\oplus } ( b * a )\, $ in the graded direct sum 
$\, \mathcal Z\, =\, ( A *_{\mathbb C} B ) \,\widehat{\oplus }\, ( A *_{\mathbb C } B )\, $ with grading given by the flip of the two summands. Let $\, u\in A\, $ (resp. $\, v\in B\, $) be a unitary and define an ungraded $C^*$-automorphism $\, {\lambda }_u\, $ of $\, \mathcal Z\, $ by letting $\, {\lambda }_u\, $ be equal to the identity in the first direct summand as above and defining 
$$ {\lambda }_u ( b_1 * a_1 * b_2 * a_2 * \cdots * b_n * a_n )\quad =\quad 
 u^* * b_1 * ua_1u^* * b_2 * ua_2u^* *\cdots * b_n * u a_n \>  $$
on the second direct summand. Define the ungraded $C^*$-automorphism $\, {\rho }_u\, $ to be equal to the identity on the second summand and given by
$$ {\rho }_u ( a_1 * b_1 * a_2 * b_2 *\cdots * a_n * b_n )\quad =\quad 
a_1u * b_1 * u^*a_2u * b_2 * \cdots * u^*a_nu * b_n * u^*  \>  $$
on the first summand. Similarly one defines $\, {\lambda }_v\, $ and $\, {\rho }_v\, $ exchanging the roles of the first and second summand. Then the left and right action of the unitary group of $\, A\, $ (resp. 
$\, B\, $) on $\, \mathfrak Z\, $ extends to a left and right action on $\, \mathcal Z\, $ putting 
$$ l_u ( x )\> =\>  u\, {\lambda }_u ( x )\quad ,\quad r_u ( x )\> =\> {\rho }_u ( x )\, u \> . $$
Also note that with respect to the grading and superinvolution the following identities hold
$$ {\lambda }_u ( x^* )\> =\> {\rho }_{u^*} ( x )^*\quad ,\quad {\lambda }_u ( \alpha ( x ) )\> =\> 
\alpha ( {\rho }_{u^*} ( x ) ) \> . $$
Then consider the super operator systems 
$\, \overline{\mathfrak Z}_A = M_2 ( A )\, {\otimes }_{h^*}\, B\, $ in case of the $A$-action, and similarly 
$\,\overline{\mathfrak Z}_B = A\, {\otimes }_{h^*}\, M_2 ( B )\, $ for the $B$-action with representations and enveloping graded $C^*$-algebras $\, \overline{\mathcal Z}_A\, $ and 
$\,\overline{\mathcal Z}_B\, $ as above. For a unitary $\, u\in A\, $ and $\, U = u \oplus 1\, $ the map 
$\, ad_U\, $ can be extended to a graded $C^*$-automorphism of the enveloping $C^*$-algebra  
$\, \overline{\mathcal Z}_A\, $. Consider in particular the pair of selfadjoint unitaries 
$$ E\> =\>  \begin{pmatrix} 1 & 0 \\ 0 & -1 \end{pmatrix} \otimes {\bf 1}_B\in M_2 ( \mathbb C )  $$
and  
$$ F\> =\>  \begin{pmatrix} 0 & 1 \\ 1 & 0 \end{pmatrix} \otimes {\bf 1}_B\in 
M_2 ( \mathbb C )\> . $$ 
Let $\, {\mathcal J}_E\, $ denote the graded $C^*$-ideal of $\, \overline{\mathcal Z}_A\, $ generated by the relations 
$$ \bigl\{\, {\lambda }_E ( x ) - x\, ,\, {\rho }_E ( x ) - x\,\bigm\vert\, x\in \overline{\mathfrak Z}_A\,\bigr\} $$
and check that it is invariant under $\, ad_U\, $, so that this automorphism drops to the quotient which also admits a compatible matrix decomposition by the strong action of $\, C^* ( E )\, $ induced from 
$C^*$-multiplication (which restricts to the image super operator system 
$\, \overline{\mathfrak Z}_{A , E}\, $). Then there is a natural (bijective) complete contraction
$$ \overline{\mathfrak Z}_{A , E} \, \twoheadrightarrow\, M_2 ( \mathfrak Z ) $$
induced by dividing out the relations
$$ \bigl\{\, {\lambda }_F ( x ) - x\, ,\, {\rho }_F ( x ) - x \,\bigm\vert\, x\in \overline{\mathfrak Z}_{A , E}\,\bigr\} $$
(check that $\, {\lambda }_F\, $ and $\, {\rho }_F\, $ pass to the quotient algebra since $\, F\, $ anticommutes with $\, E\, $) and commuting with the respective $F$-matrix decompositions and the maps $\, ad_U\, $, which renders a unital graded complete contraction to $\, M_2 ( {\mathfrak X}_0 )\, $, also compatible with $F$-matrix decompositions and $\, ad_U\, $. Put 
$\, {\overline{\mathfrak X}}_{2 , U} = {\overline{\mathfrak Z}}_{A , E}\, $ and $\, \overline{\mathfrak X} = 
\mathfrak Z\, $ and note that $\, ( 1 - e ) {\overline{\mathfrak X}}_{2 , U} ( 1 - e )\, $ is completely isometric to $\, \overline{\mathfrak X}\, $. The condition that the element 
$\, F_{U^n} = ad_{U^n} (F\cdot {\bf 1} )\, $ is an even selfadjoint  unitary in the injective envelope of 
$\, {\overline{\mathfrak X}}_{2 , U}\, $ for each $\, n\, $ is easily deduced from the fact that it defines an even selfadjoint unitary in the enveloping $C^*$-algebra 
$\, {\overline{\mathcal Z}}_{A , E} = {\overline{\mathcal Z}}_A / {\mathcal J}_E\, $ which is isomorphic with the graded direct sum of two copies of the amalgamated free product 
$$ \left( M_2 ( A )\, *_{C^* ( E )}\, ( B \oplus B )Ê\right)\> \widehat{\oplusÊ}\> 
\left( M_2 ( A )\, *_{C^* ( E )}\, ( B \oplus B ) \right) $$
and Proposition 15.10 of \cite{Pa}. Similarly one checks that $\, F\cdot F_{U^n} = F\, F_{U^n}\, $ for all 
$\, n\in\mathbb Z\, $. Put $\, {\mathcal Z}_{A , E ,0} = M_2 ( A )\, *_{C^* ( E )}\, ( B\oplus B )\, $ and 
$\, {\mathcal Z}_0 = A *_{\mathbb C} B\, $, and let 
$\, \overline{\mathcal V}_{2 , U ,l}\, $ and $\, \overline{\mathcal V}_{2 , U , r}\, $ be the subspaces linearly generated by ordered products $\, x * y\, $ resp. $\, y * x\, $ with $\, x\in M_2 ( A )\, ,\, y\in B\oplus B\, $.  
The last postulate is first checked on the corresponding subspaces for the natural $C^*$-surjections 
$$ {\pi }_{l , r} : {\overline{\mathcal Z}}_{A , E} \> \twoheadrightarrow M_2 ( {\mathcal Z}_0 )  $$
by projection onto the left or right summand in the graded direct sum as above.
One checks that the map $\, ad_{U^n}\, $ naturally drops to the diagonal of $\, M_2 ( {\mathcal Z}_0 )\, $ where it is given by $\, ad_{u^n}\, $ in the left upper corner and the identity map in the lower right corner. Then one defines $\, ad_{U^n}\, $ to be equal to (the projection onto the left or right direct summand of) the map $\, l_{u^n}\, $ as above in the upper right corner and equal to (the corresponding projection of) 
$\, r_{u^{-n}}\, $ in the lower left corner of $\, M_2 ( {\mathcal Z}_0 )\, $ and verifies that this map extends the (tautological) image of $\, ad_{U^n}\, $ on $\, M_2 ( \overline{\mathcal V}_{l , r} )\, $ such that the extension is completely isometric when restricted to each corner. Since such an extension is necessarily unique in the injective envelope, and since the subspace given by the linear span of products of elements of 
$\, M_2 ( \overline{\mathcal V}_{l , r} )\, $ with elements from the set $\, \{ F_{U^n} \}\, $ in the injective envelope can be realized taking the $C^*$-product in any  $C^*$-representation of 
$\, {\mathcal Z}_0\, $ one gets that also this chosen extension of $\, ad_{U^n}\, $ restricted to the corresponding subspaces of 
$\, M_2 ( {\mathcal Z}_0 )\, $ is unique with this property and must be equal to the corresponding extension in the injective envelope. 
The composite $\, ad_{U^n}\circ {\pi }_{l , r}\, $ is equal to 
$\, {\pi }_{l , r}\circ ad_{U^n}\, $ on the chosen subspaces which is an easy consequence of the fact that 
$\, ad_U\, $ is a left resp. right $C^* ( F_{U^n} )-C^*( F_{U^{n+1}} )$-module map for 
$\, M_2 ( \overline{\mathcal V}_l )\, $ resp. $\, M_2 ( \overline{\mathcal V}_r )\, $ which follows since the Haagerup tensor product $\, A\, {\otimes }_h\, B\, $ is a strong left $A$- and a strong right $B$-module,
so it is completely contractive and a $C^* ( F_{U^n} )$-module map for each $\, n\in\mathbb Z\, $ completing the argument that the $A$-action on $\, {\mathfrak X}_0\, $ is admissible. The case of $B$-action follows by symmetry. Since both actions are strongly completely contractive from the Theorem above and commute with each other, they combine to a strongly completely contractive action of $\, A\, {\otimes }_{max}\, B\, $. Putting 
$\, {\pi }_0 ( x ) = x\cdot {\bf 1}_{{\mathfrak X}_0}\, $ renders a unital complete contraction 
$$ {\pi }_0 :\> A\, {\otimes }_{max}\, B\> \longrightarrow\> {\mathfrak X}_0 $$
which must be a complete isometry by the assumption that $\, {\mathfrak X}_0\, $ is maximal.
\par\bigskip\noindent
It is of course impossible that the operator system associated with the symmetrized Haagerup tensor product itself should be completely isometric with the maximal $C^*$-tensor product except in trivial cases. However the following Theorem gives a partial result in this direction. Let $\, A\, $ be a 
$C^*$-algebra and consider the super operator system given by the (ordinary) Haagerup tensor product 
$\,  A\, {\otimes }_h\, A\, $ with involution defined on elementary tensors by 
$$ x\otimes y\quad\mapsto\quad y^*Ê\otimes x^* \> . $$
Let the unitary group $\, {\mathcal U}_A\, $ of $\, A\, $ act diagonally from the left by 
$$ l_u ( x\otimes y )\> =\> ux \otimes uy \> . $$ 
Then the left action commutes with the associated right action, so that $\, A\, {\otimes }_h\, A\, $ becomes a (weak) ${\mathcal U}_A$-bimodule in the sense defined above. Moreover the action restricts to the subspace $\, \mathfrak Y\, $ generated by elements of the form 
$\, \{ u \otimes u\,\vert\, u\in {\mathcal U}_A \}\, $ which again is a super operator system.
Embedding it into the amalgamated free product $C^*$-algebra 
$\, A\, *_{\mathbb C}\, A\, $ with 
${\mathbb Z}_2$-grading exchanging the two natural copies of $\, A\, $ denoted 
$\, j ( A )\, $ and $\, \overline j ( A )\, $ consider its quotient by the graded $C^*$-ideal generated by the odd elements 
$$ \bigl\{ j ( v ) \overline j ( v ) - \overline j ( v ) j ( v )\,\vert\, v\in {\mathcal U}_A \bigr\}\> .  $$
Let $\, {\mathfrak X}_0\, $ be the image of $\,\mathfrak Y\, $ in the quotient which is an operator system with weakly completely contractive adjoint diagonal action of $\, {\mathcal U}_A\, $. The argument used in the proof of the following Theorem then also shows that any operator system which is smaller than 
$\, {\mathfrak X}_0\, $ and admits a weakly completely contractive action of $\, {\mathcal U}_A\, $ compatible with the corresponding diagonal action on $\, \mathfrak Y\, $ is a strong ${\mathcal U}_A$-bimodule. To obtain a specific example we use a slightly different but related construction. Namely, consider the amalgamated free product $C^*$-algebra 
$\,  \mathcal Z = M_2 ( A ) *_{\mathbb C} M_2 ( A )\, $ with grading automorphism 
given by exchanging the two natural copies $\, j ( M_2 ( A ) )\, $ and $\, \overline j ( M_2 ( A ) )\, $ of 
$\, M_2 ( A )\, $ in the manner above. Let $\, V\in M_2 ( A )\, $ denote a unitary which is of the form 
$$ \begin{pmatrix} v & 0 \\ 0 & w \end{pmatrix} $$
with $\, v , w\in {\mathcal U}_A\, $. Divide $\, \mathcal Z\, $ by the graded $C^*$-ideal generated by the images of the odd elements
$$ \{ j ( V ) \overline j ( V ) - \overline j ( V ) j ( V ) \} \> . $$
Putting $\, E_j = j ( E )\, ,\, E_{\overline j} = \overline j ( E )\, ,\, 
F_j = j ( F )\, ,\, F_{\overline j} = \overline j ( F )\, $ and $\,\overline F = F_j\, F_{\overline j}\, $ 
($\, E\, $ and $\, F\, $ are the selfadjoint unitaries of 
$\, M_2 ( \mathbb C )\, $ as above) the relations imply in particular that $\, E_j\, $ and 
$\, E_{\overline j}\, $ commute with each other.  Then the elements 
$\, p = e_j\, e_{\overline j} = j ( e )\, \overline j ( e )\, $ and  $\, q = ( 1 -  e_j )\, ( 1 - e_{\overline j} )\, $ define orthogonal projections in the quotient $\,\mathcal Y\, $. For $\, u\in {\mathcal U}_A\, $ let 
$\, U = u\oplus 1\, $ and $\, ad_U\, $ be the free product of the inner $C^*$-automorphisms 
$\, Ad\, j ( U )\, $ on $\, j ( M_2 ( A ) )\, $ and $\, Ad\, {\overline j} ( U )\, $ on $\, \overline j ( M_2 ( A ) )\, $ by conjugation with $\, U\, $ leaving $\, p\, $ and $\, q\, $ invariant. We will also have to divide by the relations 
$$  \left\{ {\overline F}_U\, p - q\, {\overline F}_U\, ,\, 
{\overline F}_U\, q - p\, {\overline F}_U\,\vert\, u\in {\mathcal U}_A \right\}  $$ 
with $\, {\overline F}_U = ad_U ( \overline F )\, $ in order that the unitary $\,\overline F\, $ intertwines the projections $\, p\, $ and $\, q\, $.
The image of the hereditary subalgebra $\, ( p + q )\, \mathcal Y\, ( p + q )\, $ with unit $\, p + q\, $ in the quotient $C^*$-algebra is denoted $\, \mathcal X\, $. 
It admits a matrix decomposition by $C^*$-multiplication with $\, p\, $ and $\, q\, $. Then consider the operator subsystem $\, {\mathfrak X}_2\, $ of $\,\mathcal X\, $ which is generated by the elements
$$ \left\{ \begin{pmatrix} F_{\overline j}\, j ( u_1 )  \overline j ( u_1 )\, F_{\overline j}  & 
\overline j ( u_2 F )\,  j ( u_2 F )  \\  
j ( v_2  F)\, \overline j ( v_2 F )  & F_j\, \overline j ( v_1 ) j ( v_1 )\, F_j \end{pmatrix} \right\} $$
with $\, u_1\, ,\, u_2\, ,\, v_1\, ,\, v_2\in {\mathcal U}_A\, $. Its lower right corner $\, {\mathfrak X}\, $ is an operator system with unit element $\, q\, $ and corresponds to the upper right corner by $C^*$-multiplication with the unitary $\,\overline F\, $ from the left.
\par\bigskip\noindent
{\bf Theorem.}\quad ("Diagonal of tensor product" theorem)\quad
With notation as above the operator system $\, \mathfrak X\, $ is a strong 
${\mathcal U}_A$-bimodule for the (left) action given by 
$$ l_u : q\, F_j\, \overline j ( v ) j ( v )\, F_j\, q\>\mapsto\> q\, F_j\, \overline j ( u v ) j ( u v )\, F_j\, q\> . $$
\par\bigskip\noindent
{\it Proof.}\quad The first aim is to exhibit that the action is well defined and (weakly) completely isometric. To this end note that $\, ad_U\, $ drops to a $C^*$-automorphism of   $\,\mathcal X\, $ which leaves invariant the upper right corner (and lower left corner) of $\, {\mathfrak X}_2\, $, so that by the composition 
$$ \mathfrak X \buildrel \sim\over\largerightarrow \overline F\, \mathfrak X 
\buildrel ad_U\over\largerightarrow 
\overline F\, \mathfrak X \buildrel\sim\over\largerightarrow \mathfrak X $$
giving the left action of $\, u\, $ one finds that $\, l_u\, $ is completely isometric. In a second step we show that the action is admissible. Then one may apply the SCES-Theorem to conclude that the action is strongly completely contractive. One easily easily finds that the natural map 
$\, {\mathfrak X}_0Ê\twoheadrightarrow \mathfrak X\, $ is completely contractive.
Let $\, \overline{\mathfrak Z}\, $ denote the super operator system 
$\, A\, {\otimes }_h\, M_2 ( \mathbb C )\, {\otimes }_h\, A\, $ with the completely antiisometric involution
$$ x\otimes c\otimes y\quad\mapsto\quad y^*\otimes c^*\otimes x^*  $$
and note that it can be embedded in the amalgamated free product $C^*$-algebra 
$\, \overline{\mathcal Z} = A *_{\mathbb C} M_2 ( \mathbb C ) *_{\mathbb C} A\, $ with grading automorphism 
given by exchanging the two natural copies $\, j ( A )\, $ and $\, \overline j ( A )\, $ of 
$\, A\, $. Using the diagonal embedding of $\, A\, $ into $\, A \oplus A \simeq A \otimes C^* ( E )\, $  
($\, E\, $ and $\, F\, $ are the selfadjoint unitaries of $\, M_2 ( \mathbb C )\, $ as above) consider the image $\, {\overline{\mathcal Z}}_E\, $ of $\,\overline{\mathcal Z}\, $ in 
the amalgamated free product 
$\, ( A \oplus A ) *_{C^* ( E )} M_2 ( \mathbb C ) *_{C^* ( E )} ( A \oplus A )\, $. 
Then $\, {\overline{\mathcal Z}}_E\, $  admits a matrix decomposition coming from the strong graded action of $\, C^* ( E )\, $ by $C^*$-multiplication. 
One has two natural projections 
$$ {\pi }_l : \, {\overline{\mathcal Z}}_E\> \twoheadrightarrow\> M_2 ( A )\, *_{C^* ( E )}\, ( A\oplus A )\> ,\quad 
{\pi }_r :\, {\overline{\mathcal Z}}_E\>\twoheadrightarrow\> ( A\oplus A )\, *_{C^* ( E )}\, M_2 ( A ) $$
where the copy $\, j ( A )\, $ always maps to the left factor, and $\,\overline j ( A )\, $ to the right factor in the amalgamated free product (in the obvious way) and a projection 
$$ \pi :\, {\overline{\mathcal Z}}_E\>\twoheadrightarrow M_2 ( A\, *_{\mathbb C}\, A ) $$
which factors over both $\, {\pi }_l\, $ and $\, {\pi }_r\, $. Let $\, {\mathcal Z}_{r , l}\, $ denote the $C^*$-algebras which are the surjective image of 
$\, {\overline{\mathcal Z}}_E\, $ under the $*$-homomorphisms $\, {\pi }_r\, $ and  $\, {\pi }_l\, $ respectively.
Then the grading automorphism on $\, {\overline{\mathcal Z}}_E\, $ drops to a 
$*$-isomorphism $\, {\mathcal Z}_r \simeq {\mathcal Z}_l\, $. For a given unitary $\, u\in A\, $ one further divides $\, {\overline{\mathcal Z}}_E\, $ by the graded $C^*$-ideal generated by the elements
$$ \left\{  j ( u )\, \overline j ( u ) - \overline j ( u ) j ( u )\> ,\quad 
j ( u )\, F -  F j ( u )\> ,\quad \overline j ( u )\, F - F\, \overline j ( u ) \right\} \> , $$
where $\, u\, $ is shorthand for the diagonal element $\, u\oplus u\in A \oplus A\, $, which quotient is denoted $\, {\mathcal X}_U\, $. Upon dividing $\, {\mathcal Z}_{l , r}\, $ by the images of the relations as above one obtains two corresponding projections 
$$ {\pi }_{U , l / r} : {\mathcal X}_U \twoheadrightarrow {\mathcal X}_{U , l / r} $$
and a $*$-isomorphism $\, {\mathcal X}_{U , l} \simeq {\mathcal X}_{U , r}\, $ induced by the grading on 
$\, {\mathcal X}_U\, $. Let $\, \overline{\mathcal X}_U\, $ be the graded direct sum of two copies of $\, {\mathcal X}_U\, $ with grading given by the product of the grading on each summand plus the flip exchanging the summands and $\, \overline{\pi }_{U , l / r} \, $ be the (nonfaithful ungraded) representations of 
$\,\overline{\mathcal X}_U\, $ given by composing $\, {\pi }_{U , l}\, $ and $\, {\pi }_{U , r}\, $ respectively with the projection onto the left resp. right direct summand in its decomposition as a graded direct sum. 
Let $\, ad_U\, $ denote the graded $C^*$-automorphism of 
$\, \overline{\mathcal X}_U\, $ defined by 
$$ p\, j ( x )\mapsto p\, j ( u x u^* )\> ,\quad q\, j ( x )\mapsto q\, j ( x )\, , $$
$$ p\, \overline j ( x ) \mapsto p\, \overline j ( u x u^* )\> ,
\quad q\, \overline j ( x )\mapsto q\,\overline j ( x )\, , $$
$$ F\,\widehat{\oplus }\, 0\> \mapsto\> 
\left( p\, j ( u )\, F\, \overline j ( u ) + q\, \overline j ( u^* )\, F\, j ( u^* ) \right)\>\widehat{\oplus } \> 0 $$
$$ 0\, \widehat{\oplus }\, F\>\mapsto\> 
0 \>\widehat{\oplus }\> 
\left( p\, \overline j ( u )\, F\, j ( u ) + q\, j ( u^* )\, F\, \overline j ( u^* )Ê\right) $$
plus the identity map of $\, C^* ( E )\, $ and check that it is well defined.
Then consider the super operator system $\, \overline{\mathfrak X}_{2 , U}\, $ which is the subspace of 
$\,\overline{\mathcal X}_U\, $ generated by elements 
$$ \left\{ \begin{pmatrix} F\, j ( u_1 ) \overline j ( u_1 )\, F  & \overline j ( u_2 )\, F\, j ( u_2 )  \\   
j ( v_2 )\, F\, \overline j ( v_2 ) & \overline j ( v_1 ) j ( v_1 ) \end{pmatrix} \widehat{\oplusÊ} 
\begin{pmatrix} F\, j ( u_1 ) \overline j ( u_1 )\, F &  j ( u_2 )\, F\, \overline j ( u_2 ) \\ 
\overline j ( v_2 )\, F\, j ( v_2 )  &  \overline j ( v_1 ) j ( v_1 ) \end{pmatrix} \right\}   $$
and let  $\,\overline{\mathfrak X}\, $ be its lower right corner. One defines the subsystems 
$\, {\Delta }_{2 , U}\, $ and $\, {\nabla }_{2 , U}\, $ as in the definition of an admissible action above.
Let $\, r_F\, $ be the order two linear map which exchanges the columns on both sides of the graded direct sum, i.e. the image of an element as above is sent to  
$$ \begin{pmatrix} F\, j ( u_2 ) \overline j ( u_2 ) F &  \overline j ( u_1 )\, F\, j ( u_1 ) \\  
j ( v_1 )\, F\, \overline j ( v_1 ) & \overline j ( v_2 )  j ( v_2 ) \end{pmatrix}
\widehat{\oplus } \begin{pmatrix} 
F\, j ( u_2 ) \overline j ( u_2 )\, F & j ( u_1 )\, F\, \overline j ( u_1 ) \\ 
\overline j ( v_1 )\, F\, j ( v_1 ) & \overline j ( v_2 ) j ( v_2 ) \end{pmatrix}  \> . $$
One checks that it is well defined.
If $\, {\rho }_F\, $ is the composition of $\, r_F\, $ with $C^*$-multiplication by the even selfadjoint unitary $\, F = F\, \widehat{\oplus }\, F\, $ from the right, it is unital, and equal to the free product of the $C^*$-endomorphisms 
$$ j ( x )\>\mapsto\> F\, j ( x )\, F\> ,\quad \overline j ( x )\>\mapsto\> \overline j ( x ) $$
on the right side of the graded direct sum, whereas on the left side it is given by the free product of the $C^*$-endomorphisms 
$$ j ( x )\>\mapsto\> \overline j ( x )\> ,\quad \overline j ( x )\>\mapsto\> F\, j ( x )\, F $$
if restricted to the off-diagonal subsystem $\, {\nabla }_{2 , U}\, $, while on the diagonal subsystem 
$\, {\Delta }_{2 , U}\, $ it is the free product of the $C^*$-endomorphisms 
$$ j ( x )\>\mapsto\> F\,\overline j ( x )\, F\> ,\quad \overline j ( x )\>\mapsto\> j ( x )\> . $$
hence $\, {\rho }_F\, $ is completely isometric (though not graded) if restricted to either of the subsystems 
$\, {\nabla }_{2 , U}\, $ or $\, {\Delta }_{2 , U}\, $. Then also $\, r_F\, $ is completely isometric on these subsystems. Similarly, let $\, l_F\, $ be the map exchanging the rows on both sides, and 
$\, {\lambda }_F\, $ the unital map with is the composite of $\, l_F\, $ with $C^*$-multiplication by 
$\, F\, $ from the left. Then it is equal to the map
$$ j ( x )\>\mapsto\> j ( x )\> ,\quad \overline j ( x )\>\mapsto\> F\,\overline j ( x )\, F $$
on the left summand of the graded direct sum, whereas on the right side it is given by 
$$ j ( x )\>\mapsto\> \overline j ( x )\> ,\quad \overline j ( x )\>\mapsto\> F\, j ( x )\, F $$
if restricted to the diagonal subsystem $\, {\Delta }_{2 , U}\, $, and 
$$ j ( x )\>\mapsto\> F\,\overline j ( x )\, F\> ,\quad \overline j ( x )\>\mapsto\> j ( x ) $$
for an element in $\, {\nabla }_{2 , U}\, $. 
Put $\, \overline{\mathfrak Y}_{2 , U} = {\nabla }_{2 , U} + 
\overline{\mathfrak X}\, $, and check that $\, {\Delta}_{2 , U}\, $ is contained in the 
$C^* ( F )$-subbimodule generated by $\, \overline{\mathfrak Y}_{2 , U}\cap {\Delta }_{2 , U}\, $ in the injective envelope of the latter subspace. It is readily seen that $\, ad_U\, $ restricts to a completely isometric $*$-linear bijection of $\, \overline{\mathfrak Y}_{2 , U}\, $ which leaves the subspace 
$\,\overline{\mathfrak X}\, $ elementwise fixed. One also verifies the conditions that the elements 
$\, F\cdot F_{U^n} = F\, F_{U^n}\, $ are evenly graded unitaries in the injective envelope (of either of the subspaces).  Further one has two unital complete contractions 
$$ {\nu }_l : {\overline{\mathfrak X}}_{2 , U}\> \twoheadrightarrow\> {\overline{\mathcal V}}_{2 , U , l}\> ,\quad {\nu }_r : {\overline{\mathfrak X}}_{2 , U}\>\twoheadrightarrow\> {\overline{\mathcal V}}_{2 , U , r} $$
onto unital operator spaces $\, {\overline{\mathcal V}}_{2 , U , l / r}\, $ which are defined to be the images of $\, {\overline{\mathfrak X}}_{2 , U}\, $ under $\, {\overline{\pi }}_{2 , U , l / r}\, $ respectively. The natural weak 
$F$-matrix decomposition on the images which is compatible with the maps $\, {\nu }_{l / r}\, $ is checked to have the additional property that the right $F$-action on $\, {\overline{\mathcal V}}_{2 , U , r}\, $ is implemented by $C^*$-multiplication with $\, F\, $ from the right, while for 
$\, {\overline{\mathcal V}}_{2 , U , l}\, $ the left action of $\, F\, $ coincides with $C^*$-multiplication (from the left). The map $\, ad_U\, $ is checked to drop to a corresponding map when restricted to the intersection of the images of $\, {\overline{\mathfrak Y}}_{2 , U}\, $ in both quotients which is completely isometric on the images of the diagonal / off-diagonal subsystems. In case of the diagonal subsystem this map is simply given by conjugation with the image of the unitary element $\, j ( U ) \overline j ( U )\, $, 
but in case of the off-diagonal subsystem one has to make a distinction between the left and right case. 
In the first instance one takes the composition of the free product of the identity map on 
$\, \overline j ( A \oplus A )\, $ and the map $\, Ad\, j ( U )\, $ on $\, j ( M_2 ( A ) )\, $, composed with conjugation by the unitary $\,  {\overline j} ( U )\, $ of the whole $C^*$-algebra (leaving the off-diagonal subsystem invariant), whereas in the second instance the situation is reversed (mirrored).  Finally one checks that restricted to the image of an enlarged corner $\, {\Delta }^{i j}_{2 , U}\, $ the maps 
$\, {\nu }_{l / r}\, $ factor over the surjection (up to complete isometry) 
$\, {\overline{\mathcal X}}_{2 , U} \twoheadrightarrow M_2 ( \overline{\mathcal X}_U )\, $ corresponding to $\, \pi : {\overline{\mathcal Z}}_E \twoheadrightarrow M_2 ( A\, *_{\mathbb C}\, A )\, $, so the relations 
$$ \{ j ( u )\, F - F\, j ( u )\, ,\, \overline j ( u )\, F - F\, \overline j ( u ) \}  $$
automatically drop out and $\, {\overline{\mathcal X}}_U\, $ is simply graded direct sum of two copies of the image of $\, A\, *_{\mathbb C}\, A\, $ divided by the $C^*$-ideal generated by the single commutator 
$$ \{ j ( u )\, \overline j ( u ) - \overline j ( u )\, j ( u ) \} \> . $$
Then one easily verifies that the norm of an element $\, x\in M_n ( {\Delta }^{i j}_{2 , U} )\, $ satisfies 
$$ \Vert x {\Vert }_n\> =\> \max\, \{ \Vert {\nu }_l  ( x ) {\Vert }_n\> ,\> \Vert {\nu }_r  ( x ) {\Vert }_n \} 
\> . $$
To complete our argument we need to connect this construction to the previous one by a $*$-linear map 
$$ p_U : {\overline{\mathfrak X}}_{2 , U}\> \twoheadrightarrow\> M_2 ( \mathfrak X ) $$
which is compatible with the $F$-matrix decomposition and the map $\, ad_U\, $ and completely contractive restricted to an enlarged corner. But then it is clear that the image of the subspace of 
$\, A\, *_{\mathbb C}\, A\, $ modulo the $*$-ideal generated by the single relation above, generated by elements 
$$ \{ \overline j ( v )\, j ( v )\, \vert\, v\in {\mathcal U}_A \} $$
surjects onto the corresponding subspace in the free product divided by all relations of the form 
$$ \{ j ( v )\, \overline j ( v ) - \overline j ( v )\, j ( v )\, \vert\, v\in {\mathcal U}_A \} \> , $$
hence, a forteriori, onto $\, \mathfrak X\, $, while on the other hand the map to  
$\, M_2 ( \overline{\mathcal X}_U )\, $ is completely isometric restricted to an enlarged corner
so its image is contained in (a copy of) $\, M_2 ( \overline{\mathfrak X} )\, $.
The rest of the statement is immediate. One concludes that the diagonal action of $\, {\mathcal U}_A\, $ on $\,\mathfrak X\, $ is a strongly completely contractive, and is unitarily implemented in the injective envelope\qed
\par\bigskip\noindent
{\it Remark.}\quad The Theorem implies that $\, \mathfrak X\, $ in fact carries the structure of a unital 
$C^*$-algebra $\, U_A\, $ which is linearly generated by the unitary group of $\, A\, $ 
(modulo $\, \{ \pm 1 \}\, $). Suppose that 
$\, A\, $ admits a (unital) augmentation $*$-homomorphism $\, \varepsilon : A \twoheadrightarrow \mathbb C\, $.
Then the free product of the identity map of $\, j ( M_2 ( A ) )\, $ and $\, \varepsilon\, $ on 
$\, \overline j ( M_2 ( A ) )\, $ is compatible with the relations for $\,\mathcal X\, $ which amounts to a 
completely positive surjection of $\, U_A\, $ onto $\, A\, $ which is a $*$-homomorphism since it induces a group isomorphism $\, {\mathcal U}_A / \{ \pm 1 \}\,\rightarrow\, {\mathcal U}_A\, $. 
In any case the assignment $\, A \rightsquigarrow U_A\, $ is clearly functorial with respect to 
$*$-homomorphisms. 
\par\bigskip\noindent 
The definition of the maximal $C^*$-tensor product has a natural extension to super-$C^*$-algebras  which we record in the following definition.
\par\bigskip\noindent
{\bf Definition 3.}\quad (i)\quad Let $\, A\, ,\, B\, $ be two (unital) super-$C^*$-algebras. The 
{\it maximal tensor product}, denoted $\, A\, {\otimes }_{max}\, B\, $ is defined to be the super-$C^*$-algebra which is the completion of the involutive algebra generated by the algebraic tensor product of 
$\, A\, $ and $\, B\, $ with product type involution in the universal representation given by taking the direct sum (up to equivalence) over all joint (unital) commuting $*$-representations of $\, A\, $ and $\, B\, $.
\par\smallskip\noindent
(ii)\quad Let $\, {\mathfrak V}\, $ and $\, {\mathfrak W}\, $ be super operator spaces with dual super operator spaces  $\, {\mathfrak V}^*\, $ and $\, {\mathfrak W}^*\, $ respectively.  One has a natural $*$-embedding of the algebraic tensor product $\, {\mathfrak V} \otimes {\mathfrak W}\, $ with product type involution into the dual of the symmetrized Haagerup tensor product  
$\, ( {\mathfrak V}^*\, {\otimes }_{h^*}\, {\mathfrak W}^* )^*\, $, which clearly defines a cross matrix norm on the algebraic tensor product (since it is smaller than the projective, but larger than the injective tensor norm). Its completion will be denoted 
$\, {\mathfrak V}\, {\otimes }_{h_*}\, {\mathfrak W}\, $ and called the {\it dual symmetrized Haagerup tensor product}.
\par\bigskip\noindent
Note that by commutativity of the symmetrized Haagerup tensor product, also the dual symmetrized Haagerup tensor product is commutative. From this one easily infers that it is the largest commutative tensor norm on the algebraic tensor product which is smaller than the Haagerup tensor norm. In particular it must be larger than the maximal super-$C^*$-tensor norm, which is also commutative 
(and associative) and smaller than the Haagerup tensor norm. \par\bigskip\noindent
Let $\, C^* ( G )\, $ denote the full $C^*$-algebra of the discrete group $\, G\, $ and 
$\, {\Delta }_n ( G )\subseteq C^* ( G )\, \widehat\otimes\,\cdots\widehat\otimes\, C^* ( G )\, $ the super operator subspace of the $n$-fold projective tensor product of $\, C^* ( G )\, $ with itself which is the closure of the linear span of diagonal elements $\,\{ g\otimes\cdots\otimes g\,\vert\, g\in G \}\, $.  From the augmentation homomorphism $\, e : C^* ( G ) \twoheadrightarrow \mathbb C\, $ one gets a natural complete contraction $\, {\Delta }_{n+1} ( G ) \twoheadrightarrow {\Delta }_n ( G )\, $ for each $\, n\in \mathbb N\, $, so that the operator matrix norms imposed on the complex group ring 
$\, \mathbb C\, G\, $ by the diagonal embedding into 
the $n$-fold projective tensor product yields an increasing sequence of matrix norms which must necessarily converge to a limit norm since it is constant for the generators $\, \{ g \}\, $. The closure with respect to this limit norm will be denoted by $\, \Delta ( G )\, $. It is easy to see from the construction that this super operator space admits a comultiplication map 
$\, \delta : \Delta ( G ) \rightarrow \Delta ( G ) \widehat{\otimesÊ}\Delta ( G )\, $, which is involutive with respect to the product type involution on the projective tensor product (which coincides on the diagonal with the involution given by $\, x \otimes y\mapsto y^*\otimes x^*\, $), and a counit for this comultiplication given by the composition 
$$ \Delta ( G ) \longrightarrow C^* ( G ) \buildrel e\over\longrightarrow \mathbb C \> . $$ 
We let $\, {\overline A}_G = \Delta ( G )^*\, $ denote the dual operator space of 
$\, \Delta ( G )\, $ which admits a dual completely antiisometric involution (compare with the final example of section 1). The whole construction can also be carried through with the Haagerup tensor product replacing the projective tensor product. We denote the corresponding increasing sequence of super operator spaces by $\, \{ {\Delta }_n^h ( G ) \}\, $ converging up to the limit $\, {\Delta }^h ( G )\, $, and write $\, {\overline A}_G^h\, $ for the dual super operator space $\, {\Delta }^h ( G )^*\, $. 
One has the following result.
\par\bigskip\noindent
{\bf Theorem 2.}\quad The super operator space $\, {\overline A}_G\, $ admits the additional structure of a unital and counital, commutative Hopf super $C^*$-algebra with multiplication (comultiplication) induced by comultiplication (multiplication) on $\, \Delta ( G )\, $. Here comultiplication of 
$\, {\overline A}_G\, $ means a completely contractive $*$-homomorphism $\,\Delta\, $ of 
$\, {\overline A}_G\, $ into the mapping space of completely bounded maps 
$\, B_G = \mathcal C\mathcal B ( \Delta ( G )\, ,\, {\overline A}_G )\, $ which again is a super $C^*$-algebra naturally containing the spatial tensor product 
$\, {\overline A}_G \check\otimes {\overline A}_G\, $. Then the term counital means that the dualization of the inclusions 
$\, 1_l\, ,\, 1_r : \Delta ( G ) \rightarrow \Delta ( G )\widehat{\otimes }\Delta ( G )\, $ with 
$\, 1_l ( x ) = 1_G\otimes x\, ,\, 1_r ( x ) = x\otimes 1_G\, $ satisfy the relations 
$$  1_l^*\circ\Delta = id_{{\overline A}_G} = 1_r^*\circ \Delta  \> . $$ 
The super operator space $\, {\overline A}^h_G\, $ admits the additional structure of a unital commutative super $C^*$-algebra with multiplication induced by comultiplication of $\, {\Delta }^h ( G )\, $. One has a natural completely contractive $*$-homomorphism $\, {\overline A}^h_G \rightarrow {\overline A}_G\, $.
Moreover there is an involutive and contractive homomorphism from the involutive Banach algebra 
$\, {\mathcal A}_G = C^* ( G )^*\, $ which is the dual Banach space of $\, C^* ( G )\, $ (with involution given by $\, {\phi }^* ( x ) = \overline{\phi ( x^* )}\, $ and multiplication induced by the comultiplication of 
$\, C^* ( G )\, $) into $\, {\overline A}_G\, $ which factors over $\, {\overline A}^h_G\, $. 
\par\bigskip\noindent
{\it Proof.}\quad We first outline the multiplication on $\, {\overline A}_G\, $ which amounts to the specification of a completely contractive $*$-linear map 
$$ {\overline A}_G\, {\otimes }_h\, {\overline A}_G \buildrel m\over\longrightarrow {\overline A}_G $$
(for the natural involution $\, x\otimes y \mapsto y^*\otimes x^*\, $ on the Haagerup tensor product) satisfying $\, m\circ ( m \otimes id ) = m\circ ( id\otimes m )\, $. The comultiplication map $\,\delta\, $ of 
$\, \Delta ( G )\, $ can be composed with the complete contraction 
$\, \Delta ( G ) \widehat{\otimes }Ê\Delta ( G ) \twoheadrightarrow \Delta ( G )\, {\otimes }_h\, 
\Delta ( G )\, $ to yield a complete isometry 
$$ \Delta ( G ) \longrightarrow \Delta ( G )\, {\otimes }_h \Delta ( G ) $$
which, if dualized, gives the multiplication map of $\, {\overline A}_G\, $ by the composition 
$$ {\overline A}_G\, {\otimes }_h\, {\overline A}_G \longrightarrow 
( \Delta ( G )\, {\otimes }_h\,\Delta ( G ))^*
\buildrel {\delta }^*\over\longrightarrow {\overline A}_G $$
the first map being completely isometric by selfduality of the Haagerup tensor product (c.f. Theorem 9.4.7 of \cite{E-R}). It is easily checked that $\, m\, $ is commutative and associative, secondly that the element of $\, {\overline A}_G\, $ determined by $\, e\, $ is a left and right unit for $\, m\, $, and that 
$\, ( x\circ y )^* = y^*\circ x^* = x^*\circ y^*\, $ (the last identity follows from commutativity). Thus 
$\, {\overline A}_G\, $ is a unital commutative abstract super $C^*$-algebra. To construct comultiplication one notes on approximating 
$\, \Delta ( G )\, $ by the super operator spaces $\, \{ {\Delta }_n ( G ) \}\, $ and using commutativity and associativity of the projective tensor product that one gets an associative "multiplication" map 
$$ \Delta ( G ) \widehat{\otimes } \Delta ( G ) \longrightarrow \Delta ( G ) $$
induced by the multiplication of $\, C^* ( G )\, $. Dualizing this map gives the desired comultiplication 
$$ {\overline A}_G \buildrel \Delta\over\longrightarrow B_G \simeq 
( \Delta ( G )\widehat{\otimes } \Delta ( G ) )^* $$ 
(Corollary 7.1.5 and Proposition 8.1.2 of \cite{E-R}). Using the same arguments as above one finds that $\, B_G\, $ is a commutative super $C^*$-algebra and that $\, \Delta\, $ and $\, 1_l^*\, ,\, 1_r^*\, $ are 
$*$-homomorphisms. Turning to $\, {\overline A}_G^h\, $ the situation is much the same. The multiplication map is obtained from the comultiplication 
$$ {\Delta }^h ( G ) \buildrel {\delta }^h\over\largerightarrow {\Delta }^h ( G ) {\otimes }_h {\Delta }^h ( G ) 
$$
which is defined just as before. The unit is again obtained from the augmentation map 
$$ {\Delta }^h ( G ) \longrightarrow C^* ( G ) \buildrel e\over\longrightarrow \mathbb C $$
and all stated properties can be checked by dualization of the corresponding properties for 
$\, {\Delta }^h ( G )\, $ and $\, \Delta ( G )\, $ (resp. $\, C^* ( G )\, $). This also renders the contractive $*$-homomorphisms $\, {\mathcal A}_G \longrightarrow {\overline A}^h_G \longrightarrow {\overline A}_G\, $ \qed 
\par\bigskip\noindent
{\it Remark.}\quad It is conceivable that there is also some kind of Hopf structure for 
$\, {\overline A}^h_G\, $. The problem is however that the Haagerup tensor product is not commutative (and the symmetrized Haagerup tensor product which gives the same metric on the diagonal is not associative). 
Nevertheless it may be that one has some kind of comultiplication homomorphism on either 
$\, {\overline A}^h_G\, $ itself or some super $C^*$-algebra inbetween $\, {\overline A}^h_G\, $ and 
$\, {\overline A}_G\, $ possibly with range into the maximal super $C^*$-tensor product (resp. some enveloping super $C^*$-algebra thereof). We write $\, A^h_G\, $ for the closure of the image of 
$\, {\mathcal A}_G\, $ in 
$\, {\overline A}^h_G\, $, and $\, A_G\, $ for the closure of its image in $\, {\overline A}_G\, $, both of which are super $C^*$-algebras.
The situation is very similar to the dual setting where one is dealing with different kinds of group algebras, with $\, {\mathcal A}_G\, $ taking the part of 
the Banach algebra $\, {\mathit l}^1 ( G )\, $, then $\, A^h_G\, $ corresponds in some sense to the full group $C^*$-algebra $\, C^* ( G )\, $, and $\, A_G\, $ vaguely to the reduced group $C^*$-algebra 
$\, C^*_r ( G )\, $.
\par\smallskip\noindent
That the completely antiisometric involution on the abstract operator algebra $\, {\overline A}_G\, $ as above is antimultiplicative is easy to check in this example, but in fact is automatic. Namely, a unital operator algebra with a given antilinear and completely antiisosmetric involution is in particular a super operator system, so its injective envelope is a $\mathbb Z / 2\mathbb Z$-graded $C^*$-algebra (it is the same as the injective envelope of the enveloping graded operator system), so the linear grading automorphism extends by rigidity to an order two bijective unital completely isometric map which is necessarily a 
$C^*$-automorphism by Proposition 1. Then the superinvolution which is the product of the adjoint map and the grading is antimultiplicative on the injective envelope. On the other hand one knows that the embedding of the operator algebra into its injective envelope is a homomorphism, so that the superinvolution must be antimultiplicative. This is a special feature of involutive operator algebras and does not hold for Banach algebras with an isometric (antilinear) involution. 
\par\bigskip\noindent
{\it Example.}\quad Let $\, G\, $ be the finite cyclic group $\, \mathbb Z / n\mathbb Z\, $. By Pontrjagin duality the dual abelian group $\, \widehat G\, $ is again isomorphic with $\, G\, $ so that 
$\, C^* ( G ) \simeq C ( \widehat G ) \simeq C ( G )\, $ and $\, {\mathcal A}_G \simeq {\mathit l}^1 ( G )\, $. As an algebra 
$\, {\mathit l}^1 ( G )\, $ is isomorphic with $\, C^* ( G )\, $ the latter of which is the unique commutative $C^*$-algebra generated by $n$ orthogonal projections. These projections are of norm one already in 
$\, {\mathit l}^1 ( G )\, $ so that if one is to obtain a contractive representation of $\, {\mathit l}^1 ( G )\, $ into $\, \mathcal B ( \mathcal H )\, $ the image is necessarily a $C^*$-algebra, and if one assumes the representation to be injective there is only one such choice up to $C^*$-isomorphism, which is 
$\, C^* ( G )\, $. So in this case one must have that both $\, A_G\, $ and $\, A^h_G\, $ are isomorphic with 
$\, C^* ( G )\, $. Note however that the dual involution on $\, {\mathcal A}_G\, $ 
(and $\, A_G\, ,\, A^h_G\, $) does not coincide with the usual involution of $\, {\mathit l}^1 ( G )\, $ but involves an additional grading on these commutative algebras induced by (linear extension of the assignment) sending a group element $\, g\, $ to $\, g^{-1}\, $. So even in this very simple case one gets (nontrivially) graded $C^*$-algebras. In the general case it cannot be expected that either $\, A_G\, $ or $\, A^h_G\, $ should be $C^*$-algebras, only operator algebras with a completely antiisometric involution. 
\par\smallskip\noindent
As another example of (noncommutative) super $C^*$-algebras, consider the enveloping super 
$C^*$-algebra of a superunitary representation of a given topological (or locally compact) group. 
A superunitary is an element $\, U\, $ of $\, \mathcal B ( \widehat{\mathcal H} )\, $ satisfying 
$\, UU^* = U^*U = 1\, $. Since these elements form a group under multiplication it makes sense to talk of a superunitary representation meaning a group homomorphism into the group of superunitaries which is continuous for a suitable topology. It should be noted that the norm of a superunitary is always larger than one and is equal to one if and only if the element is a (evenly graded) unitary. Therefore since the norm function of a superunitary has no general bound from above there is no such thing like a universal 
group super $C^*$-algebra, covering all possible superunitary representations of $\, G\, $. 
\par


\bigskip\noindent
\begin{thebibliography}{5}
  \bibitem{B-N}
        D. Blecher, M. Neal, \emph{Metric characterizations of isometries and of unital operator spaces
        and systems}, Math. arXiv (2008)
   \bigskip
   \bibitem{Bl-Ne}
         D. Blecher, M. Neal, \emph{Metric characterizations II}, Math arXiv (2012)     
   \bigskip
   \bibitem{E-R}
        E. Effros, Z.-J. Ruan, \emph{Operator Spaces}, London Mathematical Sciences Monographs - 
        New Series ¡ 23, Oxford University Press, 2000.
   \bigskip
   \bibitem{H-L-S}
        R. Haag, J. Lopuszanski, M. Sohnius, \emph{All possible generators of supersymmetries for the 
        S-matrix}, Nuclear Physics \textbf{B88} (1975), 257--274.
   \bigskip
   \bibitem{Pa}
        V. Paulsen, \emph{Completely Bounded Maps and Operator Algebras}, Cambridge Studies in        
        Advanced Mathematics 78, Cambridge University Press, 2002.   
\end{thebibliography}
\end{document}